%




\input amstex.tex

\magnification=\magstep1
\hsize=5.5truein
\vsize=9truein
\hoffset=0.5truein
\parindent=10pt
\newdimen\nagykoz
\newdimen\kiskoz
\nagykoz=7pt
\kiskoz=2pt
\parskip=\nagykoz
\baselineskip=12.7pt


\loadeufm \loadmsam \loadmsbm

\font\vastag=cmssbx10
\font\drot=cmssdc10
\font\vekony=cmss10
\font\vekonydolt=cmssi10
\font\cimbetu=cmssbx10 scaled \magstep1
\font\szerzobetu=cmss10

\font\scVIII=cmcsc8
\font\rmVIII=cmr8
\font\itVIII=cmti8
\font\bfVIII=cmbx8
\font\ttVIII=cmtt8

\def\cim#1{{\centerline{\cimbetu#1}}}
\def\szerzo#1{{\vskip0.3truein\centerline{\szerzobetu#1}}}
\def\alcim#1{{\medskip\centerline{\vastag#1}}}
\def\tetel#1#2{{{\drot#1}{\it\szukebb~#2\tagabb}}}
\long\def\biz#1#2{{{\vekony#1} #2}}
\def\kiemel#1{{\vekonydolt#1\/}}
\long\def\absztrakt#1#2{{\vskip0.4truein{\vekony#1} #2\vskip0.5truein}}
\def\szukebb{\parskip=\kiskoz}
\def\tagabb{\parskip=\nagykoz}
\def\vonal{{\vrule height 0.2pt depth 0.2pt width 0.5truein}}

\def\CC{{\Bbb C}}

\def\ds{\displaystyle}
\def\ts{\textstyle}

\def\bmfd{{Banach manifold}}
\def\bsmfd{{Banach submanifold}}
\def\hmfd{{Hilbert manifold}}
\def\hsmfd{{Hilbert submanifold}}

\def\bvbdl{{Banach vector bundle}}

\def\Ker{\hbox{\rm Ker}}
\def\re{\hbox{\rm Re}}
\def\im{\hbox{\rm Im}}

\def\Oa{{\Omega}}
\def\oa{{\omega}}
\def\UU{{\frak U}}
\def\VV{{\frak V}}
\def\WW{{\frak W}}
\def\BB{{\frak B}}
\def\GL{\text{\rm GL}}
\def\diag{\text{\rm diag}}

\def\Hom{\text{\rm Hom}}
\def\End{\text{\rm End}}

\def\vbdl{{vector bundle}}
\def\hspc{{Hilbert space}}

\def\cts{{continuous}}

\def\idml{{infinite dimensional}}

\def\idms{{infinite dimensions}}

\def\fdml{{finite dimensional}}

\def\fdms{{finite dimensions}}

\def\pscx{{pseudoconvex}}
\def\cpscx{{coordinate pseudoconvex}}
\def\scvr{{(S)-covering}}
\def\ccvr{{(C)-covering}}
\def\lc{{locally convex}}

\def\lcas{{locally convex analytic sheaf}}
\def\lcass{{locally convex analytic sheaves}}
\def\bd{{\partial}}
\def\dist{{\hbox{\rm dist}}}

\def\homo{{homomorphism}}
\def\sbs{{Schauder basis}}

\def\cubs{{countable unconditional basis}}
\def\psh{{plurisubharmonic}}
\def\pshdom{{\psh\ domination}}

\def\st{{such that}}
\def\wrt{{with respect to}}
\def\subsub{\subset\!\subset}

\def\<{{\langle}}
\def\>{{\rangle}}

\def\RR {{\Bbb R}}

\def\AA{{\Cal A}}

\def\OO {{\Cal O}}

\def\ro{\varrho}

\def\epsz{\varepsilon}
\def\fii{\varphi}
\def\fn{func\-tion}
\def\fns{func\-tions}
\def\holo{hol\-o\-mor\-phic}

\def\ssf{(S)-sheaf}
\def\ssfs{(S)-sheaves}
\def\sres{(S)-resolution}
\def\cpair{(C)-pair}

\def\mfd{manifold}

\def\ses{short exact sequence}
\def\soes{short $1$-exact sequence}
\def\stes{short $2$-exact sequence}
\def\les{long exact sequence}
\def\loes{long $1$-exact sequence}
\def\ltes{long $2$-exact sequence}
\def\sloes{short locally $1$-exact sequence}
\def\sltes{short locally $2$-exact sequence}
\def\lloes{long locally $1$-exact sequence}
\def\lltes{long locally $2$-exact sequence}

\def\nbd{neighbor\-hood}

\def\bspc{Banach space}

\def\la{\lambda}
\def\La{\Lambda}
\def\sa{\sigma}

\def\da{{\delta}}

\def\ga{{\gamma}}
\def\aa{{\alpha}}
\def\ba{{\beta}}
\def\da{{\delta}}

\def\za{\zeta}

\def\Prop{Proposition}
\def\p#1{{\Prop~#1}}

\def\Th{{Theorem}}
\def\th{{theorem}}

\def\t#1{{\Th~#1}}

{\phantom.}
\vskip0.5truein
\cim{ANALYTIC COHOMOLOGY IN A BANACH SPACE}
\szerzo{Imre Patyi\plainfootnote{${}^1$}{\rmVIII Supported in part by a Research Initiation Grant
	from Georgia State University.}}
\absztrakt{ABSTRACT.}{Let $X$ be a \bspc\ with a \cubs\ (e.g., $X=\ell_2$ Hilbert space),
	$\Oa\subset X$
	\pscx\ open.
	 We show that the sheaf cohomology groups $H^q(\Oa,S)$ vanish for $q\ge1$
	if $S$ is a member of a fairly inclusive class of sheaves of $\OO$-modules
	over $\Oa$.
	 In particular, we have the above vanishing if $S=I$ or $S=(\OO^E)^0$,
	where $I$ is the ideal sheaf of a split complex Banach submanifold $M$
	of $\Oa$, $E\to M$ is a locally trivial \holo\ \bvbdl, and $(\OO^E)^0$ is
	the zero extension to $\Oa$ of the sheaf $\OO^E$ of germs of \holo\
	sections of $E\to M$.  
	 Some applications are also given.

	 MSC 2000: 32L05, (32L10, 32L20, 32Q28, 46G20)
	
	 Key words: analytic cohomology, pseudoconvex domains, 
	            holomorphic Banach vector bundles,
		    complex Banach manifolds.


\hfill{\it Kedves Zoli \"ocs\'emnek, sz\"ulet\'esnapj\'ara.\plainfootnote{${}^2$}
{\rmVIII To my dear Younger Brother on his birthday.}}
}

\def\sA{{1}}
\def\sB{{2}}
\def\sC{{3}}
\def\sD{{4}}
\def\sE{{5}}
\def\sF{{6}}
\def\sG{{7}}
\def\sH{{8}}
\def\sI{{9}}
\def\sJ{{10}}
\def\sK{{11}}
\def\sL{{12}}


\def\tAA{{\sA.1}}
\def\tAB{{\sA.2}}
\def\tAC{{\sA.3}}
\def\tAD{{\sA.4}}

\def\eAA{{\sA.1}}
\def\eAB{{\sA.2}}
\def\eAC{{\sA.3}}
\def\eAD{{\sA.4}}

\def\tBA{{\sB.1}}
\def\tBB{{\sB.2}}

\def\eBA{{\sB.1}}
\def\eBB{{\sB.2}}

\def\tCA{{\sC.1}}
\def\tCB{{\sC.2}}
\def\tCC{{\sC.3}}
\def\tCD{{\sC.4}}
\def\tCE{{\sC.5}}

\def\tDA{{\sD.1}}
\def\tDB{{\sD.2}}
\def\tDC{{\sD.3}}
\def\tDD{{\sD.4}}

\def\eDA{{\sD.1}}
\def\eDB{{\sD.2}}
\def\eDC{{\sD.3}}
\def\eDD{{\sD.4}}
\def\eDE{{\sD.5}}

\def\tEA{{\sE.1}}
\def\tEB{{\sE.2}}
\def\tEC{{\sE.3}}
\def\tED{{\sE.4}}
\def\tEE{{\sE.5}}
\def\tEF{{\sE.6}}
\def\tEG{{\sE.7}}

\def\eEA{{\sE.1}}
\def\eEB{{\sE.2}}
\def\eEC{{\sE.3}}
\def\eED{{\sE.4}}
\def\eEE{{\sE.5}}
\def\eEF{{\sE.6}}

\def\tFA{{\sF.1}}
\def\tFB{{\sF.2}}
\def\tFC{{\sF.3}}

\def\eFA{{\sF.1}}
\def\eFB{{\sF.2}}
\def\eFC{{\sF.3}}
\def\eFD{{\sF.4}}
\def\eFE{{\sF.5}}
\def\eFF{{\sF.6}}
\def\eFG{{\sF.7}}
\def\eFH{{\sF.8}}
\def\eFI{{\sF.9}}
\def\eFJ{{\sF.10}}
\def\eFK{{\sF.11}}
\def\eFL{{\sF.12}}
\def\eFM{{\sF.13}}
\def\eFN{{\sF.14}}
\def\eFO{{\sF.15}}

\def\tGA{{\sG.1}}
\def\tGB{{\sG.2}}
\def\tGC{{\sG.3}}

\def\eGA{{\sG.1}}

\def\tHA{{\sH.1}}

\def\tHA{{\sH.1}}

\def\eJA{\sJ.1}

\def\tKA{{\sK.1}}
\def\tKB{{\sK.2}}
\def\tKC{{\sK.3}}
\def\tKD{{\sK.4}}
\def\tKE{{\sK.5}}

\def\rD{D}
\def\rG{G}
\def\rLtA{Lt1}
\def\rLtB{Lt2}
\def\rLA{L1}
\def\rLB{L2}
\def\rLC{L3}
\def\rLD{L4}
\def\rN{N}
\def\rPA{P1}
\def\rPB{P2}
\def\rPC{P3}
\def\rPD{P4}
\def\rPsA{Ps1}
\def\rPsB{Ps2}

\alcim{\sA. INTRODUCTION.}

	 It was around $55$ years ago (ca.~$1950$) when Karl Stein defined the
	notion of Stein manifolds, and Cartan, Oka, and Serre proved two fundamental,
	and long classical, theorems about Stein manifolds (and Stein spaces) called
	ever since \Th{}s~A and B.
	 One way to express (a substantial part of) \t{}B is to say that
	over $X=\CC^n$, $n\ge1$, the sheaf cohomology groups $H^q(X,S)$ vanish
	for all $q\ge1$ if $S\to X$ is a coherent analytic sheaf.
	 Sheaf cohomology vanishing theorems hold the key to many global results
	about complex manifolds, especially to those that can be solved first
	locally and then globally by patching the local solutions to global solutions.

	 In this paper we look at a class of sheaves, called sheaves of type~(S) or
	\ssfs, over suitable complex \bspc{}s that in a way mimics the class
	of coherent analytic sheaves in \fdms, and for which we can prove vanishing
	in \idms.
	 While the class of \ssfs\ is far from being as perfect as the class of
	coherent analytic sheaves in \fdms, it is, arguably, about the best class for which
	vanishing in \idms\ can be proved with current technology, and it also
	contains the most immediately geometrically relevant analytic sheaves,
	e.g., the sheaves of germs of \holo\ sections of \holo\ \bvbdl{}s, and
	ideal sheaves of split complex \bsmfd{}s.
	 The definition of the class of \ssfs\ is fairly long, so it is given
	in its own section \S\,\sD.

         Following [\rLB] by Lempert we say that \kiemel{\pshdom} holds in a
        complex \bmfd\ $\Oa$ if for every $u:\Oa\to\RR$ locally upper bounded
        there is a $\psi:\Oa\to\RR$ \cts\ and \psh\ \st\ $u(x)<\psi(x)$ for all $x\in\Oa$.

\tetel{\t\tAA.}{{\rm(Lempert, [\rLB])}
        If $X$ is a \bspc\ with a \cubs, and\/ $\Oa\subset X$ is \pscx\ open,
	then \pshdom\ holds in\/ $\Oa$.
}

	 Here we prove \Th{}s~\tAB, \tAC, and \tAD\ below.

\tetel{\t\tAB.}{Let $X$ be a \bspc\ with a \sbs, $\Oa\subset X$ \pscx\ open,
	$M\subset\Oa$ a split complex \bsmfd\ of\/ $\Oa$, $Z$ a \bspc,
	$I^Z\to\Oa$ the sheaf of germs of \holo\ \fns\/ $\Oa\to Z$ that vanish
	on $M$, $E\to M$ a locally trivial \holo\ \bvbdl, $\OO^E\to M$ the sheaf
	of germs of \holo\ sections of $E\to M$, and\/ $(\OO^E)^0\to\Oa$ the zero
	extension to\/ $\Oa$ of the sheaf\/ $\OO^E$.
	 If \pshdom\ holds in every \pscx\ open subset of\/ $\Oa\times X$,
	then the following hold.
\vskip0pt
	{\rm(a)} The sheaf $I^Z\to\Oa$ is an \ssf.
\vskip0pt
	{\rm(b)} The sheaf
	$(\OO^E)^0\to\Oa$ is an \ssf.
}

	 The main \th\ of this paper is \t\tAC\ below.

\tetel{\t\tAC.}{Let $X,\Oa$ be as in \t\tAB, and $S\to\Oa$ an \ssf.
	Then we have the following.
\vskip0pt
	{\rm(a)} The sheaf cohomology groups $H^q(\Oa,S)$ vanish for all $q\ge1$.
\vskip0pt
	{\rm(b)} There is a \ses\
$$
	0\to K\to\OO^{Z_1}\to S\to0
\tag\eAA
 $$
 	of \lcass\ over\/ $\Oa$,
	where $Z_1$ is a \bspc, and $K$ is an \ssf, \st\
	over any \pscx\ open subset $U$ of\/ $\Oa$ and for any \bspc\ $Z$
	the image of\/ $(\eAA)$ under the functor\/ $\Hom(\OO^Z,-)$ satisfies that
$$
	0\to\Hom(\OO^Z,K)\to\Hom(\OO^Z,\OO^{Z_1})\to\Hom(\OO^Z,S)\to0
\tag\eAB
 $$
 	is exact over $U$ both on the level of germs and on the level of global
	sections.
\vskip0pt
	{\rm(c)} There is a long exact sequence
$$
	\ldots\to\OO^{Z_n}\to\OO^{Z_{n-1}}\to\ldots\to\OO^{Z_1}\to S\to0
\tag\eAC
 $$
 	of \lcass\ over\/ $\Oa$,
	where $Z_n$, $n\ge1$, is a \bspc, \st\
	over any \pscx\ open subset $U$ of\/ $\Oa$ and for any \bspc\ $Z$
	the image of\/ $(\eAC)$ under the functor\/ $\Hom(\OO^Z,-)$
	satisfies that
$$\eqalign{
	\ldots\to\Hom(\OO^Z,&\OO^{Z_n})\to\Hom(\OO^Z,\OO^{Z_{n-1}})\to\cr
	\ldots
	&\to\Hom(\OO^Z,\OO^{Z_1})\to\Hom(\OO^Z,S)\to0\cr
}
\tag\eAD
 $$
 	is exact over $U$ both on the level of germs and on the level of
	global sections.
\vskip0pt
	{\rm(d)} If `(S$\,'$)-sheaves' make up any class of \lcass\
	over \pscx\ open subsets of\/ $\Oa$, and for any
	(S$\,'$)-sheaf $S$ parts\/ {\rm(a)} and\/ {\rm(b)} above holds with `\ssf'
	replaced by `(S$\,'$)-sheaf,' then any (S$\,'$)-sheaf $S\to\Oa$ is in fact
	an \ssf\ $S\to\Oa$.
}

	 \t\tAD\ below is a geometric corollary of \t\tAC\ above.

\tetel{\t\tAD.}{With the notation and hypotheses of \t\tAB\ the following hold.
\vskip0pt
	{\rm(a)} The sheaf $I^Z$ is acyclic over\/ $\Oa$. 
\vskip0pt
	{\rm(b)} Any \holo\ \fn\ $f:M\to Z$ can be extended to 
	a \holo\ \fn\ $\tilde f:\Oa\to Z$
	with $\tilde f(x)=f(x)$ for $x\in M$.
\vskip0pt
	{\rm(c)} The sheaf $(\OO^E)^0$ is acyclic over\/ $\Oa$, and thus
        the sheaf $\OO^E$ is acyclic over $M$.
\vskip0pt
	{\rm(d)} There is a \bspc\ $Z_1$ and a \holo\ \fn\ $T:\Oa\to Z_1^*$ into the
        dual \bspc\ $Z_1^*$ of $Z_1$ \st\ as point sets $M=\{x\in\Oa:T(x)=0\}$.
\vskip0pt
	{\rm(e)} For any open $U\subset\Oa$ with $M\subset U$, there is a \pscx\ open
        $\oa\subset\Oa$ with $M\subset\oa\subset U$.
\vskip0pt
	{\rm(f)} There is a \holo\ \nbd\ retraction $r:\oa\to M$, where $\oa$ is \pscx\
        open with $M\subset\oa\subset\Oa$ and $r$ is \holo\ with $r(x)=x$ for
        $x\in M$.
}

	 The proof of \t\tAC\ follows broadly the classical proof of \t{B}
	 by Cartan, Oka, and Serre, with exhaustions coming from [\rLA, \rLB]
	 and [\rPD], patching and dimension shifting from [\rPC], and 
	 amalgamation of syzygies from [\rG], [\rLtA], and [\rPD].
	  For background see [\rLA--\rLD, \rPA].

\alcim{\sB. EXHAUSTION.}

	 This section describes a way to exhaust a \pscx\ open subset $\Oa$ of a \bspc\ $X$
	that is convenient for proving vanishing results for sheaf cohomology over $\Oa$.
	 We follow here [\rLD, \S\,2].

	 We say that a \fn\ $\aa$, call their set $\AA'$, is an
	\kiemel{admissible radius \fn} on $\Oa$ if $\aa:\Oa\to(0,1)$ is \cts\
	and $\aa(x)<\dist(x,X\setminus\Oa)$ for $x\in\Oa$.
	 We say that a \fn\ $\aa$, call their set $\AA$, is an
	\kiemel{admissible Hartogs radius \fn} on $\Oa$ if $\aa\in\AA'$
	and $-\log\aa$ is \psh\ on $\Oa$.
	 Call $\AA$ \kiemel{cofinal} in $\AA'$ if for each $\aa\in\AA'$
	there is a $\ba\in\AA$ with $\ba(x)<\aa(x)$ for $x\in\Oa$.

\tetel{\p\tBA.}{Plurisubharmonic domination holds in\/ $\Oa$ if and only if $\AA$
	is cofinal in $\AA'$.
}

\biz{Proof.}{Write $\aa=e^{-u}\in\AA'$ and $\ba=e^{-\psi}\in\AA$.
	As \pshdom\ holds on $\Oa$ for $u$ \cts\ if and only if for
	$u$ locally upper bounded, the proof of \p\tBA\ is complete.
}

	 Put $B_X(x_0,r)=\{x\in X:\|x-x_0|<r\}$ for a ball in a \bspc\ $X$,
	where $x_0\in X$, and $0<r\le\infty$.

	 It will be often useful to look at coverings by balls
	$B_X(x,\aa(x))$, $x\in\Oa$, $\aa\in\AA'$, and shrink their
	radii to obtain a finer covering by balls $B_X(x,\ba(x))$,
	$x\in\Oa$, $\ba\in\AA$.

	 Let $e_n$, $n\ge1$, be a \sbs\ in the \bspc\ $(X,\|\cdot\|)$.
	 One can change the norm $\|\cdot\|$ to an equivalent
	norm so that
$
	\big\|\sum_{i=m}^nx_ie_i\big\|\le\big\|\sum_{i=M}^Nx_ie_i\big\|
 $
 	for $0\le M\le m\le n\le N\le\infty$, $x_i\in\CC$.
	 Introduce the projections $\pi_N:X\to X$,
	$\pi_N\sum_{i=1}^\infty x_ie_i=\sum_{i=1}^Nx_ie_i$, 
	$x_i\in\CC$, $\pi_0=0$, $\pi_\infty=1$, $\ro_N=1-\pi_N$, and define for
	$\aa\in\AA$ and $N\ge0$ integer the sets
$$\eqalign{
	D_N\<\aa\>&=\{\xi\in\Oa\cap\pi_N X:(N+1)\aa(\xi)>1\},\cr
	\Oa_N\<\aa\>&=\{x\in\pi_N^{-1}D_N\<\aa\>:\|\ro_N x\|<\aa(\pi_N x)\},\cr
	D^N\<\aa\>&=\pi_{N+1}X\cap\Oa_N\<\aa\>,\cr
	\Oa^N\<\aa\>&=\{x\in\pi_{N+1}^{-1}D^N\<\aa\>:\|\ro_{N+1}x\|<\aa(\pi_Nx)\},\cr
	\BB(\aa)&=\{B_X(x,\aa(x)):x\in\Oa\},\cr
	\BB_N(\aa)&=\{B_X(x,\aa(x)):x\in\Oa_N\<\aa\>\}.\cr
}
\tag\eBA
 $$

	 These $\Oa_N\<\aa\>$ are \pscx\ open in $\Oa$, and they will serve
	to exhaust $\Oa$ as $N=0,1,2,\ldots$ varies.

\tetel{\p\tBB.}{{\rm(Lempert)} Let $\aa\in\AA$, and suppose that \pshdom\ holds in\/ $\Oa$.
\vskip0pt  
	{\rm(a)} There is an $\aa'\in\AA$, $\aa'<\aa$, with\/
	$\Oa_n\<\aa'\>\subset\Oa_N\<\aa\>$ for all $N\ge n$.
	 So any $x_0\in\Oa$ has a \nbd\ contained in all but
	finitely many $\Oa_N\<\aa\>$.
\vskip0pt  
	{\rm(b)} There are $\ba,\ga\in\AA$, $\ga<\ba<\aa$, so that
	for all $N$ and $x\in\Oa_N\<\ga\>$
$$
	B_X(x,\ga(x))\subset\Oa_N\<\ba\>\cap\pi_N^{-1}B_X(\pi_N x,\ba(x))
	\subset B_X(x,\aa(x)).
\tag\eBB
 $$
\vskip0pt  
	{\rm(c)} If\/ $8\aa\in\AA$, $Y\subset X$ is a \fdml\ complex affine
	subspace, then $Y\cap\overline{\Oa_N\<\aa\>}$ is \psh{}ally convex
	in $Y\cap\Oa$.
\vskip0pt
	{\rm(d)} We have that\/ $\overline{\Oa_N\<\aa\>}\subset\overline{\Oa^N\<\aa\>}$.
	If\/ $4\aa\in\AA$, then\/ $\Oa^N\<\aa\>\subset\Oa_N\<2\aa\>$.
\vskip0pt
	{\rm(e)} There is a $\ba\in\AA$, $\ba<\aa$, with
	$\Oa_N\<\ba\>\subset\Oa_N\<\aa\>\cap\Oa_{N+1}\<\aa\>$ for $N\ge0$.
\vskip0pt
	{\rm(f)} There is an $\aa'\in\AA$, $\aa'<\aa$, \st\
	the covering $\BB_N(\aa)|\Oa_N\<\aa'\>$ has a finite basic refinement
	for all $N\ge0$.
}

\biz{Proof.}{For (a) and (b) see [\rLD, Prop.\,2.1], and [\rLC, Prop.\,4.3],
	for (c) [\rLC, Prop.\,4.3], for (d) [\rLC, Prop.\,4.4], for (e) [\rLD, Prop.\,2.3],
	and for (f) see [\rPC, Prop.\,3.2(c)].  
	 The proof of \p\tBB\ is complete.
	 (Remark for the record that (f) was not explicitly formulated by Lempert.)
}

	The meaning of \p\tBB(b) is that certain refinement maps exist
	between certain open coverings, while (cd) are useful for Runge
	type approximation, and (ef) for exhaustion.

\alcim{\sC. MODEL SHEAVES AND THEIR HOMOMORPHISMS.}

	 In this section we look at the simplest kinds of sheaves of
	$\OO$-modules, their topology on their spaces of sections, and
	their continuous homomorphisms.

	 A \kiemel{complex \bmfd} $\Oa$ modelled on a \bspc\ $X$ is a paracompact
	Hausdorff space with an atlas of bi\holo{}ally related charts onto open
	subsets of $X$.
	 Many of the complex analytic properties of $\Oa$
	can be studied by looking at
	the sheaves $\OO^Z\to\Oa$ of germs of \holo\ \fns\ $\Oa\to Z$, where
	$Z$ is any \bspc.
	 We call any such sheaf $\OO^Z$ a model sheaf over $\Oa$.

	 The vector space $\OO(U,Z)=\OO^Z(U)$ of global sections of $\OO^Z$
	over any open $U\subset\Oa$ carries a natural complete locally convex
	vector topology induced by the family of seminorms $p_K$, where $K$ runs
	through all compact subsets of $U$, defined by $p_K(f)=\sup_{x\in K}\|f(x)\|_Z$
	for $f\in\OO(U,Z)$.
	 As the point evaluations $\OO(U,Z)\ni f\mapsto f(x)\in Z$, $x\in U$,
	are \cts\ linear functionals in this topology, the
	space $\OO(U,Z)$ is indeed Hausdorff.
	 If $\Oa$ is \fdml, then the resulting locally convex spaces $\OO(U,Z)$
	are in fact Fr\'echet spaces.
	 If $\Oa$ is \idml, then $\OO(U,Z)$ may not be a Fr\'echet spaces.

	 We denote by $H=\Hom(\OO^{Z_1},\OO^{Z_2})$ the sheaf of $\OO$-linear
	\cts\ sheaf homomorphisms over $\Oa$ from $\OO^{Z_1}$ to $\OO^{Z_1}$,
	i.e., the sections $\tau\in H(U)$ over an open $U\subset\Oa$ are
	$\OO$-linear maps $\tau:\OO^{Z_1}|U\to\OO^{Z_2}|U$ that induce
	\cts\ linear maps $\tau:\OO(V,Z_1)\to\OO(V,Z_2)$ in the natural
	topology discussed above for all open $V\subset U$.
	 (For our purposes it is enough to consider topology on $\OO(V,Z)$
	only for arbirarily small \cpscx\ open \nbd{}s $V$ of every point of $\Oa$.)

	 Any \holo\ operator \fn\ $T\in\OO(\Oa,\Hom(Z_1,Z_2))$ induces a sheaf
	\homo\ $\tau=\dot T\in\Hom(\OO^{Z_1},\OO^{Z_2})\to\Oa$ defined by
	$\dot Tf=g$, where $f\in\OO(U,Z_1)$, $g\in\OO(U,Z_2)$, and
	$g(x)=T(x)f(x)$ for $x\in U$, $U\subset\Oa$ open.

	 We go on to show that any $\tau$ arises as $\tau=\dot T$ for a
	unique $T$ at least if $X$ is nice enough.

\tetel{\p\tCA.}{Let\/ $\Oa$ be a complex \bmfd\ modelled on a \bspc\ $X$ with a
	\sbs, $Z_1,Z_2$ \bspc{}s, and $\tau:\OO^{Z_1}\to\OO^{Z_2}$ a sheaf \homo\
	of $\OO$-module sheaves over\/ $\Oa$.
	 Suppose that $\tau$ is sequentially \cts, i.e.,
	if for each point $x_0\in\Oa$ there is a \cpscx\ open set $V$ with
	$x_0\in V\subset\Oa$ \st\ for every \cpscx\ open $U$ with $x_0\in U\subset V$, 
	$f_n,f\in\OO(U,Z_1)$, $n=1,2,3,\ldots,$ 
	and $f_n\to f$ uniformly
	on compact subsets of $U$ as $n\to\infty$, then
	$\tau f_n\to\tau f$ uniformly on compact subsets of $U$ as $n\to\infty$.
	 Then $\tau$ is of the form $\tau=\dot T$ for a unique
	$T\in\OO(\Oa,\Hom(Z_1,Z_2))$.
}

	 Note that for any such $T$ the induced sheaf \homo\ $\dot T$ is
	(sequentially) \cts, since $\|T(x)\|$ is bounded for $x$ in any compact
	subset $K$ of $\Oa$.
	 The proof of \p\tCA\ will occupy us for a while.

\tetel{\p\tCB.}{{\rm(a)} Let $X$ be a \bspc\ with a \sbs, and Schauder projections
	$\pi_n,\ro_n:X\to X$ as in\/ \S\,\sB.
	 Then for any $x\in X$ the sequence\/ $\|\ro_n(x)\|$ decreases down to zero
	as $n\to\infty$.
\vskip0pt
	{\rm(b)} If $f_n:K\to[0,\infty)$ are \cts\ \fns\ on a compact space $K$,
	and for each $x\in K$ the sequence $f_n(x)$ decreases down to zero as
	$n\to\infty$, then $f_n\to0$ uniformly on $K$.
\vskip0pt
	{\rm(c)} For any compact $K\subset X$ we have that\/ $\sup_{x\in K}\|\ro_n(x)\|\to0$
	as $n\to\infty$.
\vskip0pt
	{\rm(d)} Let $B_X(x_0,R)$ be an open ball in $X$, $Z$ a \bspc, and
	$f:B_X(x_0,R)\to Z$ \holo, or, even just \holo\ on all complex affine
	one dimensional slices of\/ $B_X(x_0,R)$.
	 If $f$ is bounded on $B_X(x_0,R)$, then $f$ is Lipschitz \cts\ on
	$B_X(x_0,r)$ for $r<R/4$.
\vskip0pt
	{\rm(e)} If $f\in\OO(B_X(0,r),Z)$ is Lipschitz \cts, then $f_n\!\in\!\OO(B_X(0,r),Z)$
	defined by $f_n(x)=f(\pi_n(x))$ tends to $f$ uniformly on compact subsets
	of $B_X(0,r)$.
}

\biz{Proof.}{(a) See the paragraph of $(\eBA)$.

	(b) This is a classical theorem of Dini.
	 Given any $\epsz>0$, for any $x_0\in K$ let $N_\epsz(x_0)$ be the
	smallest index $N\ge1$ with $f_N(x_0)<\epsz$.
	 As the inequality $f_N(x)<\epsz$ persists for $x$ in an open \nbd\ of $x_0$,
	we see that the \fn\ $x_0\mapsto N_\epsz(x_0)$ is locally upper bounded on $K$.
	 A locally upper bounded \fn\ on a compact space is in fact globally upper 
	bounded.
	 Hence there is an integer $M_\epsz\ge1$ \st\ $N_\epsz(x)\le M_\epsz$
	for all $x\in K$, i.e.,
	for all $n\ge M_\epsz$ we have that $0\le f_n(x)\le f_{M_\epsz}(x)<\epsz$
	for all $x\in K$, or, $f_n\to0$ uniformly on $K$.

	(c) By (a) part (b) is applicable to $f_n(x)=\|\ro_n(x)\|$.

	(d) By assumption there is an $M>0$ with $\|f(x)\|<M$ for $\|x\|<R$.
	 If $\|x\|<r$, $\|y\|<r$, $x\not=0$, $|\la|<3r$, then
	$\big\|x+\la\frac{y-x}{\|y-x\|}\big\|\le\|x\|+|\la|<r+3r=4r<R$.
	 Let $z^*\in Z^*$ be a \cts\ linear functional on the \bspc\ $Z$ with
	$\|z^*\|\le1$.
	 Thus the \fn\ $\fii(\la)=z^*\big(f(x+\la\frac{y-x}{\|y-x\|})-f(x)\big)$,
	$|\la|<3r$, is a numerical \holo\ \fn, and satisfies that
	$\fii(0)=0$, and $|\fii(\la)|\le2M$ for $|\la|<3r$.
	 The classical Schwarz lemma implies that
	 $|\fii(\la)|\le\frac{2M}{3r}|\la|$ for $|\la|<3r$.
	 In particular, for $\la=\|y-x\|<2r$ we have that
	$|\fii(\|y-x\|)|=|z^*(f(y)-f(x))|\le\frac{2M}{3r}\|y-x\|$.
	 Taking supremum for $z^*\in Z^*$ with $\|z^*\|\le1$ we find that
	$\|f(y)-f(x)\|\le\frac{2M}{3r}\|y-x\|$ as claimed.

	(e) Let $L$ be a Lipschitz constant for $f$, i.e.,
	$\|f(y)-f(x)\|\le L\|y-x\|$ for $\|x\|,\|y\|<r$.
	As $\|f_n(x)-f(x)\|\le\|f(\pi_n(x))-f(x)\|\le L\|\pi_n(x)-x\|\le
	L\|\ro_n(x)\|$ an application of (c) completes the proof of \p\tCB.
}

\tetel{\p\tCC.}{If $T:\Oa\to\Hom(Z_1,Z_2)$ satisfies that $T(x)z_1$ is \holo\
	in $x\in\Oa$ for each fixed $z_1\in Z_1$, then $T$ is \holo, i.e.,
	$T\in\OO(\Oa,\Hom(Z_1,Z_2))$.
}

\biz{Proof.}{This is a well-known classical statement.
	 As the desired holomorphy of $T$ is a local property, we may assume
	that $\Oa=B_X(0,1)$.
	 If $x_n\to x_0$ in $\Oa$, then $T(x_n)z_1\to T(x_0)z_1$.
	So by, say, the Banach--Steinhaus theorem $\|T(x_n)\|$ is bounded
	as $n\to\infty$, i.e., $\|T(x)\|$ is a locally bounded \fn\ of
	$x\in\Oa$.
	 
	 We need to show that $T$ is \holo\ on one dimensional complex affine
	slices of $\Oa$.
	 Fix $x_0,x_1\in X$, $\|x_1\|=1$, and look at the \fn\ 
	$T_\la=T(x_0+\la x_1)$
	of $\la$ on the open set $\La$ of the $\la$-plane where it is defined.
	 Before we can show the desired holomorphy of $T_\la$ for
	$\la\in\La$ we check that $T_\la$ is locally Lipschitz
	\cts\ on $\La$.
	 To that end let $M$ be a bound of $\|T(x)\|$ in a \nbd\ in $\Oa$ of the
	point $x_0+\la_0 x_1$ for any fixed $\la_0\in\La$, and $z_1\in Z_1$
	any vector with $\|z_1\|\le1$.
	 On applying the Schwarz lemma as in the proof of \p\tCB(d) to the
	\holo\ \fn\ $T_\la z_1$ in a small disc about $\la_0$
	we see that an estimate
	$\|T_{\la_2} z_1-T_{\la_1}z_1\|\le L|\la_2-\la_1|$
	holds, where $L$ is independent of $z_1$ and $\la_1,\la_2$.
	 Taking supremum for $z_1\in Z_1$, $\|z_1\|\le1$, we obtain that
	$\|T_{\la_2}-T_{\la_1}\|\le L|\la_2-\la_1|$, i.e.,
	$T_\la$ is indeed locally Lipschitz \cts\ for $\la\in\La$.
	 Thus the vector valued Riemann integral
	$S_\mu=\frac{1}{2\pi i}\int_{|\la-\la_0|=\epsz}\frac{T_\la}{\la-\mu}\,d\la$
	exists for $\mu$ in a small disc $D=B_{\CC}(\la_0,\epsz)$ with
	positively oriented boundary circle.
	 Clearly, $S_\mu$ is \holo\ for $\mu\in D$.
	 As $T_\la z_1$ is \holo\ in $\la\in D$ for each fixed $z_1\in Z_1$, it
	satisfies the Cauchy integral formula
	$T_\mu z_1=\frac{1}{2\pi i}\int_{|\la-\la_0|=\epsz}\frac{T_\la z_1}{\la-\mu}\,d\la=
	S_\mu z_1$.
	 Thus the operator \fn\ $T_\mu$ equals the \holo\ operator \fn\
	$S_\mu$ for $\mu\in D$.
	 So far we have seen that $T$ is locally bounded, and one dimensional
	slicewise \holo\ on $\Oa$.
	 By \p\tCB(d) our $T$ is also locally Lipschitz \cts\ on $\Oa$.
	 It is a simple classical fact that if a \fn\ $T$ is \cts\ and
	\holo\ on one dimensional complex affine slices of $\Oa$, then
	$T$ is \holo.
	 The proof of \p\tCC\ is complete.
}

\tetel{\p\tCD.}{Let $X,Z_1,Z_2$ be \bspc{}s, $B=B_X(0,1)$ a ball,
	$e_i\in X$, $\xi_i\in X^*$ \cts\ linear functionals, $i=1,\ldots,n$,
	and $\tau:\OO^{Z_1}\to\OO^{Z_2}$ a sheaf \homo\ of $\OO$-module
	sheaves, \cts\ or not.
	 If $f(0)=0$ for an $f\in\OO(B,Z_1)$, then the \fn\ $g\in\OO(U,Z_1)$
	defined in an open \nbd\ of\/ $0\in X$ by $g(x)=f(\sum_{i=1}^n\xi_i(x)e_i)$
	satisfies that $(\tau g)(0)=0$.
}

\biz{Proof.}{Look at the \fn\ $h(\za)=f(\sum_{i=1}^n\za_i e_i)$ defined and
	\holo\ for $\za=(\za_1,\ldots,\za_n)$ in an open \nbd\ of the
	origin in $\CC^n$.
	 As $h(0)=0$ we can write $h(\za)=\sum_{i=1}^n\za_i h_i(\za)$
	in a \nbd\ of the origin, where $h_i$, $i=1,\ldots,n$, is \holo\
	in a \nbd\ of the origin, either by power series expansion, or
	by looking at $h(\za)=h(\za)-h(0)=\int_{t=0}^1\frac{d}{dt}h(t\za)\,dt$.
	 Hence $g(x)=\sum_{i=1}^n\xi_i(x)g_i(x)$ for $x$ in a \nbd\ of
	the origin of $X$, where the $g_i$ are \holo.
	 Since $(\tau g)(x)=\sum_{i=1}^n\xi_i(x)(\tau g_i)(x)$, on setting
	$x=0$ we get that $(\tau g)(0)=0$.
	 The proof of \p\tCD\ is complete.
}

\tetel{\p\tCE.}{With the notation and hypotheses of \p\tCA\ if $f\in\OO(U,Z_1)$,
	$f(x_0)=0$, $x_0\in U\subset\Oa$ open, then $(\tau f)(x_0)=0$.
}

\biz{Proof.}{Without loss of generality we may assume that $x_0=0\in X$,
	$U=B_X(0,1)$ and $f$ is Lipschitz \cts\ on $U$.
	 Define $f_n\in\OO(U,Z_1)$ by $f_n(x)=f(\pi_n(x))$.
	 Since $f_n(0)=0$ \p\tCD\ implies that $(\tau f_n)(0)=0$ for $n\ge1$.
	\p\tCB(e) shows that $(\tau f_n)\to(\tau f)$ uniformly on compact
	subsets $K$ of $U$.
	 In particular, letting $K=\{0\}$ yields that
	$0=(\tau f_n)(0)\to(\tau f)(0)$ as $n\to\infty$.
	 Thus $(\tau f)(0)=0$, and the proof of \p\tCE\ is
	complete.
}

\biz{Proof of \p\tCA.}{Letting $f_{z_1}$ be various constant \fns\
	$f_{z_1}(x)=z_1\in Z_1$ we see that $\tau(f_{z_1})(x)=T(x)z_1$
	for a unique linear map $T(x):Z_1\to Z_2$.
	 If a sequence $z_1^n\to z_1$ converges in $Z_1$ in norm, then
	$f_{z_1^n}\to f_{z_1}$ uniformly on (compact subsets of) $X$,
	so $T(x)\in\Hom(Z_1,Z_2)$ is indeed a bounded linear operator,
	and $T\in\OO(\Oa,\Hom(Z_1,Z_2))$ by \p\tCC.
	 It remains to show that $\tau=\dot T$, i.e.,
	$\tau f=\dot Tf$ for $f\in\OO(U,Z_1)$, $U\subset\Oa$ open.
	 \p\tCE\ says that if $f(x_0)=0$ for an $x_0\in U$, then
	$(\tau f)(x_0)=0$ as well.
	 For a general $f$ write $f(x)=(f(x)-f(x_0))+f(x_0)$.
	 If we regard $f(x_0)$ as a constant member of $\OO(U,Z_1)$,
	then $\tau f=\tau(f-f(x_0))+\dot Tf(x_0)$, whose value at
	$x_0$ is $(\tau f)(x_0)=0+T(x_0)f(x_0)$.
	 Hence $\tau=\dot T$ and the proof of \p\tCA\ is complete.
}

	 Note that in the above proof of \p\tCA\ the sequential
	continuity of $\tau$, and topology on section spaces,
	were used only on arbitrarily small
	open neighborhoods $U$ of any point $x_0$ of $\Oa$, where
	$U$ is bi\holo\ to a ball in $X$.

\alcim{\sD. (S)-SHEAVES.}

	 In this section we define a class of analytic sheaves called
	\ssfs\ that form the major object of study in this paper, and
	we also look at some of their first properties.

	 Let $\Oa$ be a complex \bmfd\ modelled on a \bspc\ $X$.
	 We call an open $U\subset\Oa$ \kiemel{\cpscx}
	if $U$ is bi\holo\ to a \pscx\ open subset of $X$.
	An open subset $U$ of $X$ is \cpscx\ if and only if it is \pscx. 
	 We call an open covering $\UU$ of $\Oa$ an \kiemel{\scvr} if all
	intersections $\bigcap U_i$ of finitely many members $U_i$ of
	$\UU$ are \cpscx.
	 It is easy to see that any open covering $\VV$ of $\Oa$ has
	a refinement $\UU$ that is an \scvr.

	 Let $S\to\Oa$ be a sheaf of $\OO$-modules.
	 We call $S$ a \kiemel{\lcas} over $\Oa$ if $\Oa$ has an \scvr\ $\UU$ \st\
	for any $U\in\UU$, $V\subset U$ \cpscx\ open, the set of sections
	$S(V)$ carries a complete Hausdorff \lc\ topological vector space
	structure so that the $\OO$-module multiplication $\OO(V)\times S(V)\to S(V)$
	is \cts, and the restriction maps $V_2\subset V_1$ induce \cts\
	linear maps $S(V_1)\to S(V_2)$.

	 Let $S_1$ and $S_2$ be two \lcass\ over $\Oa$.
	 Let $\UU_i$ be an \scvr\ of $\Oa$ that can serve in the above definition
	for $S_i$, $i=1,2$.
	 Let $\UU$ be an \scvr\ of $\Oa$ that is a common refinement of $\UU_1$
	and $\UU_2$.
	 Let $\tau:S_1\to S_2$ be a sheaf \homo\ of $\OO$-module sheaves over $\Oa$.
	 We call $\tau$ \kiemel{\cts} over $\Oa$ and write $\tau\in\Hom(S_1,S_2)(\Oa)$
	if for each $U\in\UU$, $V\subset U$ \cpscx\ open, the map induced by $\tau$
	on sections $S_1(V)\to S_2(V)$ is \cts\ in the given topologies of $S_1$
	and $S_2$.
	 The set of all \cts\ sheaf \homo{}s $\tau$ as above form a sheaf
	$\Hom(S_1,S_2)$ over $\Oa$, whose sections over any open $U\subset\Oa$ are
	all the \cts\ sheaf \homo{}s $S_1|U\to S_2|U$ of \lcass.
	 This sheaf $\Hom(S_1,S_2)$ may or may not be a \lcas\ over $\Oa$.

	 The model sheaf $\OO^Z\to\Oa$, where $Z$ is any \bspc, is, with its natural
	topology, a \lcas\ over $\Oa$, and the sheaf of \cts\ \homo{}s
	$\Hom(\OO^{Z_1},\OO^{Z_2})$, where $Z_1,Z_2$ are \bspc{}s, is naturally
	isomorphic to the sheaf of \holo\ operator \fns\ $\OO^{\Hom(Z_1,Z_2)}$
	by \p\tCA\ if $X$ has a \sbs.
	 The main reason to look at topology on spaces of sections, and continuity
	of sheaf \homo{}s is precisely the above identification of 
	$\Hom(\OO^{Z_1},\OO^{Z_2})$ with $\OO^{\Hom(Z_1,Z_2)}$; a triviality
	if both $Z_1$ and $Z_2$ are \fdml, as in the case of classical Stein theory.

	 The only \lcass\ $S$ that will interest us in this paper are locally
	of the form $\OO^Z/K$, where $K$ is a subsheaf of $\OO^Z$ with $K(U)$
	being a closed subspace of $\OO(U,Z)$, where $U$ is a small enough
	\cpscx\ open \nbd\ of any point of $\Oa$.
	 Such sheaves $S$ are indeed \lcass\ over $\Oa$.

	 Let $X$ be a \bspc, $\Oa\subset X$ \pscx\ open, and
	$M$ a complex \bmfd\ modelled on $X$.
	 Let
$$
	0\to A\to B\to C\to0
\tag\eDA
 $$
 	be a \ses\ of $\OO$-module sheaves over $\Oa$.

	 We say that $(\eDA)$ is \kiemel{$1$-exact} or a
	\kiemel{\soes} over a \pscx\ open subset $U$ of $\Oa$ if $(\eDA)$
	is exact on the germ level at any point $x\in U$, and on the level
	of global sections over any \pscx\ open $V\subset U$.

	 Let $(\eDA)$ be a \ses\ of $\OO$-module sheaves over $M$.
	 We say that $(\eDA)$ is \kiemel{locally $1$-exact} or a
	\kiemel{\sloes} over $M$ if there is an \scvr\ $\UU$ of $M$ \st\
	for all $U\in\UU$ our $(\eDA)$ is a \soes\ over $U$ in the above sense.

	 Let $(\eDA)$ be a \ses\ of \lcass\ (and their \cts\ \homo{}s) over $\Oa$.
	 We say that $(\eDA)$ is \kiemel{$2$-exact} or a \kiemel{\stes} over
	a \pscx\ open $U\subset\Oa$ if for any \bspc\ $Z$ the image of $(\eDA)$
	under the functor $\Hom(\OO^Z,-)$ satisfies that
$$
	0\to\Hom(\OO^Z,A)\to\Hom(\OO^Z,B)\to\Hom(\OO^Z,C)\to0
\tag\eDB
 $$
 	is a \soes\ over $U$.

	 Let $(\eDA)$ be a \ses\ of \lcass\ over $M$.
	 We say that $(\eDA)$ is \kiemel{locally $2$-exact} or a \kiemel{\sltes}
	over $M$ if there is an \scvr\ $\UU$ of $M$ \st\ $(\eDA)$ is a \stes\
	over each $U\in\UU$ in the above sense.

	 Let $i=1,2$, and
$$
	\ldots\to A_n{\buildrel{\tau_n}\over\to}A_{n-1}\to\ldots\to
	A_1{\buildrel{\tau_1}\over\to}A_0\to0
\tag\eDC
 $$
	a \les\ of sheaves of $\OO$-modules over $\Oa$ for $i=1$, and a \les\
	of \lcass\ over $\Oa$ for $i=2$, $K_p=\Ker\,\tau_n$, $p\ge1$, $K_0=A_0$,
	the associated kernel sheaves, and
$$
	0\to K_{p+1}\to A_{p+1}\to K_p\to0,
\tag\eDD
 $$
 	$p\ge0$, the sequence of \ses{}s over $\Oa$ associated to $(\eDC)$.

	 If $(\eDD)$ is a short $i$-exact sequence over a \pscx\ open $U\subset\Oa$
	for all $p\ge0$, then we say that $(\eDC)$ is \kiemel{$i$-exact} or a
	\kiemel{long $i$-exact sequence} over $U$.

	 Let $i=1,2$, and $(\eDC)$ a \les\ of sheaves of $\OO$-modules over $M$
	for $i=1$, and a \les\ of \lcass\ over $M$ for $i=2$, and let $K_p$,
	and $(\eDD)$ as above.
	 If there is an \scvr\ $\UU$ of $M$ so that for all $U\in\UU$ our $(\eDD)$
	is a short $i$-exact sequence over $U$ for all $p\ge0$, then we say that
	$(\eDC)$ is \kiemel{locally $i$-exact} or a \kiemel{long locally $i$-exact
	sequence} over $M$.

	 Let $Z_n$, $n\ge1$, be a \bspc, $S\to M$ a sheaf of $\OO$-modules.
	 A \les\
$$
	\ldots\to\OO^{Z_n}\to\OO^{Z_{n-1}}\to\ldots\to\OO^{Z_1}\to S\to0
\tag\eDE
 $$
 	of sheaves of $\OO$-modules over $M$ is called a \kiemel{$1$-resolution}
	(by model sheaves) over $M$ if $(\eDE)$ is a \lloes\ over $M$.

	 Let $Z_n$, $n\ge1$, be a \bspc, $S\to M$ a \lcas.
	A \les\ $(\eDE)$ of \lcass\ over $M$ is called a \kiemel{$2$-resolution}
	or an \kiemel{\sres} (by model sheaves) over $M$ if $(\eDE)$ is a \lltes\
	over $M$.

	 We call a \lcas\ $S\to M$ a \kiemel{sheaf of type~(S)} or an \kiemel{\ssf}
	if there is an \scvr\ $\UU$ of $M$ \st\ over each $U\in\UU$ our sheaf $S|U$
	admits an \sres\ as above.

	 To deal with \sres{}s and \ssfs\ \t\tED\ and the following two theorems
	come in handy.

\tetel{\t\tDA.}{{\rm([\rPD, Thm.\,1.3])}
	 Let $X$ be a \bspc\ with a \sbs, $\Oa\subset X$ \pscx\ open,
	$E\to\Oa$ a \holo\ \bvbdl\ with a \bspc\ $Z$ for fiber type.
	 If \pshdom\ holds in\/ $\Oa$, then we have the following.
\vskip0pt
	{\rm(a)} $H^q(\Oa,\OO^Z)=0$ for $q\ge1$.
\vskip0pt
	{\rm(b)} Let $Z_1=\ell_p(Z)$, $1\le p<\infty$.
	 Then $E\oplus(\Oa\times Z_1)$ and\/ $\Oa\times Z_1$ are \holo{}ally
	isomorphic over\/ $\Oa$.
\vskip0pt
	{\rm(c)} $H^q(\Oa,\OO^E)=0$ for $q\ge1$.
\vskip0pt
	{\rm(d)} If $E$ is \cts{}ly trivial over\/ $\Oa$, then $E$ is \holo{}ally
	trivial over\/ $\Oa$.
}

\tetel{\t\tDB.}{{\rm([\rPC, Thm.\,4.3])}
	 Let $X$ be a \bspc\ with a \sbs, $\Oa\subset X$ \pscx\ open, and suppose
	that \pshdom\ holds in every \pscx\ open subset of\/ $\Oa$.
	 Let $(\eDE)$ be a\/ $1$-resolution of a sheaf $S$ of $\OO$-modules
	over\/ $\Oa$, and $K_p$, $p\ge0$, the associated sequence of
	kernel sheaves over\/ $\Oa$.
	 Then the following hold.
\vskip0pt
	{\rm(a)} $H^q(\Oa,K_p)=0$ for all $q\ge1$ and $p\ge0$.
\vskip0pt
	{\rm(b)} The sequence\/ $(\eDE)$ is exact over\/ $\Oa$ on the level of
	global sections.
}

	 Next we show that locally exact global resolutions are in fact
	globally exact.

\tetel{\t\tDC.}{Let $X$ be a \bspc\ with a \sbs, $\Oa\subset X$ \pscx\ open,
	$i=1,2$, $S\to\Oa$ a sheaf of $\OO$-modules for $i=1$, and a \lcas\
	for $i=2$.
	 Suppose that \pshdom\ holds in every \pscx\ open subset of\/ $\Oa$.
	 If\/ $(\eDE)$ is an $i$-resolution of $S$ over\/ $\Oa$, $K_p$, $p\ge0$,
	are the associated kernel sheaves, then $K_p$, $p\ge0$, are acyclic over
	any \pscx\ open $U\subset\Oa$.
	 In particular, $(\eDE)$ is a long $i$-exact sequence over\/ $\Oa$,
	and for $i=2$ the sheaves $K_p$, $p\ge0$, are \ssfs\ over\/ $\Oa$.
}

\biz{Proof.}{For $i=1$ this follows from \t\tDB.
	 For $i=2$ this follows from the case $i=1$ above applied to
	$(\eDE)$, and then to its image under the functor $\Hom(\OO^Z,-)$,
	where $Z$ is any \bspc.
	 Finally, the sheaf $K_p$ is an \ssf, since
$$
	\ldots\to\OO^{Z_n}\to\OO^{Z_{n-1}}\to\ldots\to\OO^{Z_{p+1}}\to K_p\to0
 $$
 	is an \sres\ of $K_p$, $p\ge0$, over $\Oa$,
	where all maps are as in $(\eDE)$.
	 The proof of \t\tDC\ is complete.
}

	 Our goal is to show eventually that any \ssf\ $S\to\Oa$ as in \t\tDC\
	admits in fact an \sres\ over all of $\Oa$.

\tetel{\p\tDD.}{{\rm(a)} The direct sum of two \sres{}s
$$\eqalign{
	&\ldots\to\OO^{Z'_n}\buildrel{\tau'_n}\over\to\OO^{Z'_{n-1}}\to\ldots
	\to\OO^{Z'_1}\buildrel{\tau'_1}\over\to S'\to0,\cr
	&\ldots\to\OO^{Z''_n}\buildrel{\tau''_n}\over\to\OO^{Z''_{n-1}}\to\ldots
        \to\OO^{Z''_1}\buildrel{\tau''_1}\over\to S''\to0\cr
} $$
	over the same complex \bmfd\ $M$
	is an \sres\/ $(\eDE)$ over $M$, where $Z_n=Z'_n\oplus Z''_n$,
	$\tau_n={{\tau'_n\,\,0}\choose{0\,\,\tau''_n}}$, $n\ge1$, and $S=S'\oplus S''$.
\vskip0pt
	{\rm(b)} If in an \sres\/ $(\eDE)$ over a complex \bmfd\ $M$ the map
	$\tau_1:\OO^{Z_1}\to S$ is replaced by $\tau_1A:\OO^{Z_1}\to S$, where
	$A\in\OO(M,\GL(Z_1))$, and the other maps $\tau_n$, $n\ge2$, are
	unchanged, then we get another \sres\ of $S$ over $M$.
}

\biz{Proof.}{As both parts are clear from the definitions, the proof of \p\tDD\
	is complete.
}
 	
\alcim{\sE. EXAMPLES OF (S)-SHEAVES.}

	 In this section we show that some of the simplest and most immediately
	geometrically relevant analytic sheaves are \ssfs.

	 In \idms\ no analog of the classical Oka coherence theorem seems to be
	currently available.
	 It is therefore difficult to verify whether a sheaf is an \ssf, requiring
	a case by case study, and we can do it here in only a few, but useful, cases.

	 The definition of an \ssf\ is purely local.
	 We only need to exhibit an \sres\ on small enough \cpscx\ open \nbd{}s of
	each point of the ground \bmfd.

	 As an aside, note that over a \fdml\ complex \mfd\ any coherent analytic
	sheaf is an \ssf, as it is easy to see using the basic theorems of Stein
	theory, such as the Oka coherence theorem and \t{B}.
	 Also there are many \ssfs\ that are not of a finite rank, let alone
	coherent analytic.

\tetel{\p\tEA.}{Let $M$ be a complex \bmfd, and $Z$ a \bspc.
	Then $\OO^Z\to M$ is an \ssf.
}

\biz{Proof.}{As it is easy to check that the trivial sequence
	$\ldots\to0\to0\to\ldots\to0\to\OO^Z\buildrel1\over\to\OO^Z\to0$
	is a \lltes\ over $M$, the proof of \p\tEA\ is complete.
}

	 Note that to make any use of \ssfs\ one has to have a thorough
	understanding of the model sheaves $\OO^Z$, in particular, one
	has to know that over any \pscx\ open subset of the ground \bspc\
	the sheaves $\OO^Z$ are acyclic.

\tetel{\p\tEB.}{Let $M$ be a complex \bmfd, and $E\to M$ a \holo\ \bvbdl\
	with a \bspc\ $Z$ for fiber type.
	 Then the sheaf $\OO^E$ of germs of \holo\ sections $E\to M$ is an
	\ssf\ over $M$.
}

\biz{Proof.}{Restricting $E$ to members of an \scvr\ of $M$ by coordinate balls
	$U$ over which $E$ is \holo{}ally isomorphic to $U\times Z$, we see that
	it is enough to apply \p\tEA\ to $\OO^E|U$ to conclude 
	the proof of \p\tEB.
}

\tetel{\t\tEC.}{Let $X$ be a \bspc\ with a \sbs, $\Oa\subset X$ \pscx\ open.
	 Suppose that \pshdom\ holds in every \pscx\ open subset of\/ $\Oa\times X$.
	  Let $Z$ be a \bspc, $X=X'\times X''$ a direct decomposition of \bspc{}s,
	 and $I^Z\to X$ the sheaf of germs of \holo\ \fns\ $X\to Z$ that vanish
	 on $X'$.
	  Then $I^Z$ is an \ssf\ over\/ $\Oa$.
}

	 The canonical Koszul resolution of $I^Z$ is an \sres, as we will see
	shortly.
	 In fact, this Koszul resolution served as model for the notion of \sres.

	 Let $\La^{\!Z}_p$ the \bspc\ of all \cts\ complex $p$-linear
	alternating maps $X''\to Z$ for $p\ge0$; $\La^{\!Z}_0=\La^{\!Z}_{-1}=Z$; 
	and $\OO^{\La^{\!Z}_p}\to X$ the
	sheaf of germs of \holo\ functions $X\to\La^{\!Z}_p$.
	 Let $E$ be the Euler vector 
	field on $X''$ defined by $E(x'')=x''$, 
	and $i_E$ the inner
	derivation determined by the vector field $E$, i.e., $i_E$ is the contraction
	of $p$-forms with $E$: if $f$ is a local section of $\OO^{\La^{\!Z}_p}$, then
	let $i_Ef$ be the local section of $\OO^{\La^{\!Z}_{p-1}}$ given for $p\ge1$ by
	$(i_Ef)(x',x'')(\xi''_1,\xi''_2,\ldots,\xi''_{p-1})=
	 f(x',x'')(x'',\xi''_1,\ldots,\xi''_{p-1})$, and for $p=0$ by
	$(i_Ef)(x',x'')=f(x',0)$.
	 We consider the Koszul complex
$$
	\ldots\to\OO^{\La^{\!Z}_p}\to\OO^{\La^{\!Z}_{p-1}}\to\ldots\to
	\OO^{\La^{\!Z}_1}\to I^Z\to0
\tag\eEA
 $$
 	of \lcass\ over $X$, where each map is $i_E$.
	 Let $K_p$, $p\ge0$, be the corresponding sequence of kernel sheaves:
	$K_p(U)=\{f\in\OO(U,\La^Z_p):i_Ef=0{\text{\/ on\ }}U\}$, $U\subset X$ open; $K_0=I^Z$.

\tetel{\t\tED.}{{\rm[\rPC, Thm.\,5.1]} 
	 Let $X',X'',Z$ be \bspc{}s, $\Oa\subset X=X'\times X''$ \pscx\ open,
	$I^Z$ the sheaf of germs of \holo\ \fns\ $\Oa\to Z$ that vanish on $X'$.
	 Suppose that $X$ has a \sbs, and that \pshdom\ holds in every \pscx\ open
	subset of\/ $\Oa\times X$. Then
{\vskip0pt\rm (a)}
	the Koszul complex\/ $(\eEA)$ is exact on the germ level and on the level
	of global sections over\/ $\Oa$, and
{\vskip0pt\rm (b)}
	the $K_p$ are acyclic over\/ $\Oa$: $H^q(\Oa,K_p)=0$ for all $q\ge1$ and $p\ge0$.
}

\biz{Proof of \t\tEC.}{Look at the image
$$\eqalign{
	\ldots\to&\OO^{\Hom(Z',\La^{\!Z}_p)}\to\OO^{\Hom(Z',\La^{\!Z}_{p-1})}\to\ldots\cr
	&\phantom{12345}\to\OO^{\Hom(Z',\La^{\!Z}_1)}\to\Hom(\OO^{Z'},I^Z)\to0\cr
}
\tag\eEB
 $$
 	of $(\eEA)$ under the functor $\Hom(\OO^{Z'},-)$, where $Z'$ is any \bspc.
	 Let $Z''=\Hom(Z',Z)$, and note that the \bspc{}s $\Hom(Z',\La^{\!Z}_p)$
	and $\La^{\!Z''}_p$ are canonically isomorphic for $p\ge0$, and that the
	sheaves $\Hom(\OO^{Z'},I^Z)$ and $I^{Z''}$ are canonically isomorphic over $X$.
	 Moreover, the sequence $(\eEB)$ is canonically isomorphic to the sequence
$$
	\ldots\to\OO^{\La^{\!Z''}_p}\to\OO^{\La^{\!Z''}_{p-1}}\to\ldots\to
	\OO^{\La^{\!Z''}_1}\to I^{Z''}\to0
\tag\eEC
 $$
 	over $X$, which, being just another Koszul complex, is a \loes\ over $X$
	by \t\tED, i.e., $(\eEA)$ is a \ltes\ over $X$.
	 The proof of \t\tEC\ is complete.
}

\tetel{\t\tEE.}{Let $M$ be a complex \bmfd\ modelled on a \bspc\ $X$ with a \sbs,
	and suppose that \pshdom\ holds in any \pscx\ open subset of $B_X(0,1)\times X$.
	 Let $N$ be a closed split complex \bsmfd\ of $M$, $Z$ a \bspc, and $I^Z$
	the sheaf of germs of \holo\ \fns\ $M\to Z$ that vanish on $N$.
	 Then $I^Z$ is an \ssf\ over $M$.
}

\biz{Proof.}{Let $\UU$ be an \scvr\ of $M$ so that if $U\in\UU$ meets $N$, then
	the pair $(U,U\cap N)$ is bi\holo\ to the pair 
	$(B_{X'}(0,1)\times B_{X''}(0,1),B_{X'}(0,1)\times\{0\})$, where
	$X=X'\times X''$ is a direct decomposition of \bspc{}s.
	 As $I^Z|U$ has an \sres\ by \t\tEC, the proof of \t\tEE\ is complete.
}

	 We now turn to zero extensions of analytic sheaves from split complex
	\bsmfd{}s, and show in many cases that the zero extension of an \ssf\
	over the submanifold is an \ssf\ over the ambient \mfd.

	 Let $M$ be a complex \bmfd, $N$ a split complex \bsmfd\ of $M$,
	$S\to N$ a sheaf of $\OO$-modules over the submanifold $N$.
	 The \kiemel{zero extension} or \kiemel{trivial extension}
	$S^0\to M$ of $S\to N$ is the sheaf of the canonical presheaf $S^0$
	defined over $M$ by letting $S^0(U)=S(U\cap N)$, $U\subset M$ open.
	 Note that $S(\emptyset)=0$ as usual.

\tetel{\t\tEF.}{Let $X$ be a \bspc\ with a \sbs, $X=X'\times X''$ a direct
	decomposition of \bspc{}s, $\Oa=B_{X'}(0,1)\times B_{X''}(0,1)$,
	$\Oa'=B_{X'}(0,1)\times\{0\}$, regard\/ $\Oa'$ as a split complex \bsmfd\
	of\/ $\Oa$, and let $S\to\Oa'$ be an \ssf\ with a \ltes\/ $(\eDE)$ of \lcass\
	over\/ $\Oa'$.
	 If \pshdom\ holds in any \pscx\ open subset of\/ $\Oa\times X$,
	then the zero extension $S^0\to\Oa$ of $S\to\Oa'$ has a \ltes\ of model
	sheaves over $\Oa$.
}

\biz{Proof.}{Denote the maps in $(\eDE)$ by $\tau'_n:\OO^{Z_n}\to\OO^{Z_{n-1}}$
	for $n\ge2$, and $\tau'_1:\OO^{Z_1}\to S$ over $\Oa'$.
	 \p\tCA\ gives a $T'_n\in\OO(\Oa',\Hom(Z_n,Z_{n-1}))$ with 
	 $\tau'_n=\dot T'_n$ for $n\ge2$.
	 Let $T_n(x',x'')=T'_n(x')$ be the trivial extension
	$T_n\in\OO(\Oa,\Hom(Z_n,Z_{n-1}))$ of $T'_n$ for $n\ge2$.
	 Let $i_E$ be the inner derivation 
	$i_E\in\OO(X,\Hom(\La^{\!Z_n}_p,\La^{\!Z_n}_{p-1}))$ of the Koszul complex
	$(\eEA)$ for $n\ge1$, $p\ge1$.
	 Note that
$$
	T_{n-1}T_n=0,\quad T_mi_E=i_ET_m,\quad i_Ei_E=0,
\tag\eED
 $$
 	where $n\ge3$ in the first part, $T_mi_Ef=i_ET_mf$ for $m\ge2$,
	$f\in\OO(U,\La^{\!Z_n}_p)$, $n\ge1$, $p\ge1$, $U\subset\Oa$ open,
	in the second part, since both sides have the value
	$T'_m(x')f(x',x'')(x'',\xi''_1,\ldots,\xi''_{p-1})$ at a point
	$(x',x'')\in U$.

	 Let $\tilde Z_n=\bigoplus_{p=1}^n\La^{\!Z_p}_{n-p}$ for $n\ge1$,
	and define $\tilde\tau_n:\OO^{\tilde Z_n}\to\OO^{\tilde Z_{n-1}}$
	for $n\ge2$, and $\tilde\tau_1:\OO^{\tilde Z_1}\to S^0$ over $\Oa$
	by letting $\tilde\tau_1(f_1)=\tau'_1(f_1|\Oa'\cap U)\in S(\Oa'\cap U)=
	S^0(U)$ for $f_1\in\OO(U,\tilde Z_1)$, $\tilde Z_1=Z_1$, $U\subset\Oa$
	\pscx\ open, and $\tilde\tau_n(f_1,\ldots,f_n)=(g_1,\ldots,g_{n-1})$,
	$n\ge2$, by
$$
	g_k=T_{k+1}f_{k+1}+(-1)^ki_Ef_k,
 $$
 	where $f_k\in\OO(U,\La^{\!Z_k}_{n-k})$ for $k=1,\ldots, n$, and
	$g_k\in\OO(U,\La^{\!Z_k}_{n-1-k})$ for $k=1,\ldots,n-1$.
	 Then $\tilde\tau_1\tilde\tau_2(f_1,f_2)=\tilde\tau_1(T_2f_2-i_Ef_1)=
	\tau'_1(\tau'_2(f_2|\Oa'\cap U)+0)=0$ since $(\eDE)$ is a complex,
	$\tilde\tau_{n-1}\tilde\tau_n(f_1,\ldots,f_n)=(h_1,\ldots,h_{n-2})$,
	where $n\ge3$, $h_k=T_{k+1}g_{k+1}+(-1)^ki_Eg_k=
	T_{k+1}(T_{k+2}f_{k+2}+(-1)^{k+1}i_Ef_{k+1})+(-1)^ki_E(T_{k+1}f_k)=
	0+(-1)^{k+1}(T_{k+1}i_Ef_{k+1}-i_ET_{k+1}f_{k+1})+0=0$, taking $(\eED)$
	into account.

	 Let now $g_1\in\Ker\,\tilde\tau_1$, i.e., $\tau'_1(g_1|\Oa'\cap U)=0$,
	where $g_1\in\OO(U,\tilde Z_1)$, $\tilde Z_1=\La^{\!Z_1}_0=Z_1$, $U\subset\Oa$
	\pscx\ open.
	 As $(\eDE)$ is exact on the level of global sections over $\Oa'\cap U$, there is a
	representation $g_1|\Oa'\cap U=T'_2f'_2$, where $f'_2\in\OO(\Oa'\cap U,Z_2)$.
	 \t\tED(b) gives an extension $f_2\in\OO(U,Z_2)$ with $f_2=f'_2$ on $\Oa'\cap U$.
	 Then $g_1-T_2f_2$ vanishes on $\Oa'$, so by \t\tED(a) it can be written as
	$g_1-T_2f_2=-i_Ef_1$, where $f_1\in\OO(U,\La^{\!Z_1}_1)$, i.e.,
	$g_1=\tilde\tau_2(f_1,f_2)$.

	 Let now $f\in\Ker\,\tilde\tau_n$, $n\ge2$.
	 Consider the system of equations
$$\left\{\eqalign{
	g_1&=T_2f_2-i_Ef_1=0\cr
	g_2&=T_3f_3+i_Ef_2=0\cr
	&\phantom{1}\vdots\cr
	g_{n-1}&=T_nf_n+(-1)^{n-1}i_Ef_{n-1}=0\cr
}\right.,
 $$
 	for the unknowns $f_n,f_{n-1},\ldots,f_1$.
	 We must show that all solutions are of the form
$$\left\{\eqalign{
	f_1&=T_2h_2-i_Eh_1\cr
	f_2&=T_3h_3+i_Eh_2\cr
	&\phantom{1}\vdots\cr
	f_n&=T_{n+1}h_{n+1}+(-1)^ni_Eh_n\cr
}\right..
 $$
	 We will look at the equations $g_{n-1}=0$, \dots, $g_1=0$ in this
	order and produce the representatives $h_{n+1},\ldots,h_1$ of
	$f_n,\ldots,f_1$ in this order.
	 Restricting the equation $g_{n-1}=0$ to $\Oa'\cap U$ we find that
	$T_n(f_n|\Oa'\cap U)=0$, i.e., $f_n|\Oa'\cap U=T_{n+1}h'_{n+1}$,
	where $h'_{n+1}\in\OO(\Oa'\cap U,\La^{\!Z_{n+1}}_0)$ since $(\eDE)$
	is exact on the level of global sections over $\Oa'\cap U$ by assumption
	and \t\tDC.
	 \t\tED(b) provides an extension $h_{n+1}\in\OO(U,\La^{\!Z_{n+1}}_0)$
	\st\ $h_{n+1}=h'_{n+1}$ on $\Oa'\cap U$.
	 Thus $f_n-T_{n+1}h_{n+1}$ vanishes on $\Oa'\cap U$.
	 \t\tED(a) then provides an $h_n\in\OO(U,\La^{\!Z_{n+1}}_1)$ with
	$f_n=T_{n+1}h_{n+1}+(-1)^ni_Eh_n$.
	 Looking at the equation $g_{n-1}=0$ again, we get that
	$0=g_{n-1}=(-1)^ni_E(T_nh_n-f_{n-1})$,
	where we plugged in the above form of $f_n$, and used the first two
	identities in $(\eED)$.
	 As $T_nh_n-f_{n-1}$ is in the kernel of $i_E$ over the \pscx\ open set $U$,
	another application of \t\tED(a) gives an $h_{n-1}$ with
	$f_{n-1}=T_nh_n+(-1)^{n-1}i_Eh_{n-1}$.
	 Plugging this form of $f_{n-1}$ into the equation $g_{n-2}=0$ we obtain
	the compatibility condition that $0=g_{n-2}=(-1)^{n-1}i_E(T_{n-1}h_{n-1}-f_{n-2})$,
	so as above we find an $h_{n-2}$ with $f_{n-2}=T_{n-1}h_{n-1}+(-1)^{n-2}i_Eh_{n-2}$.
	Continuing in this way we find $h_{n+1},h_n,\ldots,h_1$ one after another
	by \t\tED(a).
	 Thus any solution of $\tilde\tau_n(f)=0$ is of the form $f=\tilde\tau_{n+1}(h)$.
	 Hence our sequence
$$
	\ldots\to\OO^{\tilde Z_n}\buildrel{\tilde\tau_n}\over\to\OO^{\tilde Z_{n-1}}\to
	\ldots\to\OO^{\tilde Z_1}\buildrel{\tilde\tau_1}\over\to S^0\to0
\tag\eEE
 $$
 	is a \loes\ of \lcass\ over $\Oa$.
	 Taking the image of $(\eEE)$ under the functor $\Hom(\OO^Z,-)$, where $Z$ is
	any \bspc, we obtain a similar sequence
$$\eqalign{
        \ldots\to&\OO^{\Hom(Z,\tilde Z_n)}\buildrel{\tilde\tau_n}\over\to
	\OO^{\Hom(Z,\tilde Z_{n-1})}\to \ldots\cr
	&\phantom{12345}\to\OO^{\Hom(Z,\tilde Z_1)}
	\buildrel{\tilde\tau_1}\over\to\Hom(\OO^Z,S^0)\to0,\cr
}
\tag\eEF
 $$
 	which is a \loes\ over $\Oa$ by essentially the same reasoning as above for
	$(\eEE)$ noting that in all our arguments we could carry a linear parameter
	$z\in Z$.
	 The proof of \t\tEF\ is complete.
}

\tetel{\t\tEG.}{Let $M$ be a complex \bmfd\ modelled on a \bspc\ $X$ with a \sbs,
	and suppose that \pshdom\ holds in any \pscx\ open subset of $B_X(0,1)\times X$.
	 Let $N$ be a closed split complex \bsmfd\ of $M$, $S\to N$ and \ssf, and
	$S^0\to M$ the zero extension of $S\to N$.
	 Then $S^0$ is an \ssf\ over $M$.
}

\biz{Proof.}{This follows from \t\tEF\ in a way similar to the proof of \t\tEE.
	 The proof of \t\tEG\ is complete.
}

\alcim{\sF. AMALGAMATION OF SYZYGIES.}

	 In this section we paste together resolutions over neighboring \pscx\
	open sets.

	 Let $X$ be a \bspc, we call a pair of \pscx\ open subsets $U',U''$ of $X$
	a \kiemel{\cpair} if $U'\cup U''$ is also \pscx\ open in $X$.
	 A fairly typical example of a \cpair\ can be obtained as follows.
	 Let $\Oa\subset X$ be \pscx\ open, $f\in\OO(\Oa)$, and
	$-\infty<a'<a''<b'<b''<\infty$ constants, and define
	$U'=\{x\in\Oa:a'<\re\,f(x)<b'\}$ and
	$U''=\{x\in\Oa:a''<\re\,f(x)<b''\}$.
	 Then $U',U''$ is a \cpair\ in $X$.

\tetel{\t\tFA.}{Let $X$ be a \bspc\ with a \sbs, $U',U''\subset X$ a \cpair,
	$U=U'\cup U''$, $V=U'\cap U''$, $S\to U$ a \lcas.
	 Suppose that \pshdom\ holds in every \pscx\ open subset of\/ $U$.
\vskip0pt
	{\rm(a)} If $S|U'$ has an \sres\ over $U'$, and
	$S|U''$ has an \sres\ over $U''$, then there is a \ses\
$$
	0\to K\to\OO^{Z_1}\to S\to0
\tag\eFA
 $$
	of \lcass\ over $U$, where $Z_1$ is a \bspc, \st\ the restriction of\/
	$(\eFA)$ to $U'$ is a \stes\ over $U'$, and the resriction of\/ $(\eFA)$ to $U''$
	is a \stes\ over $U''$, and $K|U'$ has an \sres\ over $U'$, and
	$K|U''$ has an \sres\ over $U''$.
\vskip0pt
	{\rm(b)} There is an \sres\/ $(\eDE)$ over $U$.
}

\biz{Proof.}{Let

\noindent
$(\eFB)$\hfil\hbox{
$\ds
	\ldots\to\OO^{Z'_n}\buildrel{\tau'_n}\over\to\OO^{Z'_{n-1}}\to
	\ldots\to\OO^{Z'_1}\buildrel{\tau'_1}\over\to S|U'\to0,
 $}\hfil

\noindent
$(\eFC)$\hfil\hbox{
$\ds
	\ldots\to\OO^{Z''_n}\buildrel{\tau''_n}\over\to\OO^{Z''_{n-1}}\to
	\ldots \to\OO^{Z''_1}\buildrel{\tau''_1}\over\to S|U''\to0
 $}\hfil

\noindent
	be \sres{}s over $U'$ and $U''$, and

\noindent
$(\eFD)$\hfil\hbox{
$\ds\eqalign{
        &\ldots\to\OO^{\Hom(Z''_1,Z'_n)}\buildrel{\tau'_n}\over\to
	\OO^{\Hom(Z''_1,Z'_{n-1})}\to\cr
        &\phantom{12345}\ldots\to\OO^{\Hom(Z''_1,Z'_1)}\buildrel{\tau'_1}\over\to
	\Hom(\OO^{Z''_1},S)|U'\to0,\cr
 }
 $}\hfil

 \noindent
$(\eFE)$\hfil\hbox{
$\ds\eqalign{
         &\ldots\to\OO^{\Hom(Z'_1,Z''_n)}\buildrel{\tau''_n}\over\to
	 \OO^{\Hom(Z'_1,Z''_{n-1})}\to\cr
         &\phantom{12345}\ldots\to\OO^{\Hom(Z'_1,Z''_1)}\buildrel{\tau''_1}\over\to 
	 \Hom(\OO^{Z'_1},S)|U''\to0\cr
 }
 $}\hfil

\noindent
	the image of $(\eFB)$ under the functor $\Hom(\OO^{Z''_1},-)$, and the image
	of $(\eFC)$ under the functor $\Hom(\OO^{Z'_1},-)$.
	 As $(\eFD)$ is exact, by \t\tDC, on the level of global sections over $V$
	there is for $\tau''_1\in\Hom(\OO^{Z''_1},S)(V)$ an
	$A\in\OO(V,\Hom(Z''_1,Z'_1))$ \st\
	$\tau'_1\dot A=\tau''_1$ over $V$.
	 As $(\eFE)$ is exact, by \t\tDC, on the level of global sections over $V$
	there is for $\tau'_1\in\Hom(\OO^{Z'_1},S)(V)$ a
	$B\in\OO(V,\Hom(Z'_1,Z''_1))$ \st\
	$\tau''_1\dot B=\tau'_1$ over $V$.
	 (The above intertwining property is the main reason to look at the
	notion of $2$-exactness.)

	 Let the direct sum of $(\eFB)$ with the trivial \sres

\noindent
\hfil\hbox{
$\ds 
	\ldots\to0\to0\to\OO^{Z''_1}\buildrel1\over\to\OO^{Z''_1}\to0\to0
 $}\hfil

\noindent
	over $U'$ be

\noindent
$(\eFF)$\hfil\hbox{
$\ds\eqalign{
        &\ldots\to\OO^{Z'_n}\buildrel{\tau'_n}\over\to\OO^{Z'_{n-1}}\to \ldots\cr
	&\phantom{12345}\to\OO^{Z'_3}\buildrel{{\tau'_3}\choose{0}}\over\longrightarrow
	\OO^{Z'_2\oplus Z''_1}\buildrel{{\tau'_2\,\,0}\choose{0\,\,\,\,1}}\over\longrightarrow
	\OO^{Z'_1\oplus Z''_1}\buildrel{(\tau'_1,0)}\over\longrightarrow
	S|U'\to0,\cr
 }
 $}\hfil

\noindent
	which is another \sres\ by \p\tDD(a).

	 Let the direct sum of the trivial \sres

\noindent
\hfil\hbox{
$\ds
	\ldots\to0\to0\to\OO^{Z'_1}\buildrel1\over\to\OO^{Z'_1}\to0\to0
 $}\hfil

\noindent
	with $(\eFC)$ over $U''$ be

\noindent
$(\eFG)$\hfil\hbox{
$\ds\eqalign{
	&\ldots\to\OO^{Z''_n}\buildrel{\tau''_n}\over\to\OO^{Z''_{n-1}}\to\ldots\cr
	&\phantom{12345}\to\OO^{Z''_3}\buildrel{{0}\choose{\tau''_3}}\over\longrightarrow
	\OO^{Z'_1\oplus Z''_2}\buildrel{{\!\!\!1\,\,\,0}\choose{0\,\,\tau''_2}}\over\longrightarrow
	\OO^{Z'_1\oplus Z''_1}\buildrel{(0,\tau''_1)}\over\longrightarrow
	S|U''\to0,\cr
 }
 $}\hfil

\noindent
	which is another \sres\ by \p\tDD(a).

	 Let $Z_1=Z'_1\oplus Z''_1$, and consider the \holo\ operator \fn\ 
	$C_t\in\OO(V,\GL(Z_1))$ defined for $0\le t\le1$ by
$$
	C_t(x)=\bmatrix
	       1& tA(x)\cr
	       0& 1\cr
	       \endbmatrix
	       \bmatrix
	       1& 0\cr
	       -tB(x)& 1\cr
	       \endbmatrix,
 $$
	whose inverse is 
	$C_t(x)^{-1}=
	\big[{{1\ \ \ \ \ \ 0}\atop{tB(x)\,\,1}}\big]
        \big[{{1\ -tA(x)}\atop{0\ \ \ \ \ \ \ 1}}\big]$.
	 Note that $(\tau'_1,0)C_1=(0,\tau''_1)$ over $V$.
	 Look at the \holo\ \bvbdl\ $E\to U$ whose fiber type is
	$Z_1$, and whose transition \fn\ relative to the open covering
	$\{U',U''\}$ of $U$ is $C_1$.
	 This $E$ is \cts{}ly trivial over $U$ due to the homotopy $C_t$
	from $C_1$ to $C_0=1$.
	 By \t\tDA(d) our \holo\ \bvbdl\ $E$ is \holo{}ally trivial over $U$.
	 So there are \holo\ operator \fns\
	$C'\in\OO(U',\GL(Z_1))$, and
	$C''\in\OO(U'',\GL(Z_1))$ with
	$C_1(x)=C'(x)C''(x)^{-1}$ for $x\in V$.
	 Thus the \homo{}s
	$(\tau'_1,0)C'\in\Hom(\OO^{Z_1},S)(U')$,
	and $(0,\tau''_1)C''\in\Hom(\OO^{Z_1},S)(U'')$
	fit together to a \homo\
	$\tau_1\in\Hom(\OO^{Z_1},S)(U)$.
	
	 Replace $(\tau'_1,0)$ in $(\eFF)$ by $\tau_1=(\tau'_1,0)C'$ to obtain

\noindent
$(\eFH)$\hfil\hbox{
$\ds\eqalign{
        &\ldots\to\OO^{Z'_n}\buildrel{\tau'_n}\over\to\OO^{Z'_{n-1}}\to \ldots\cr
        &\phantom{12345}\to\OO^{Z'_3}\buildrel{{\tau'_3}\choose{0}}\over\longrightarrow
        \OO^{Z'_2\oplus Z''_1}\buildrel{{\tau'_2\,\,0}\choose{0\,\,\,\,1}}\over\longrightarrow
        \OO^{Z_1}\buildrel{\tau_1}\over\longrightarrow
        S|U'\to0.\cr
 }
 $}\hfil

	 Replace $(0,\tau''_1)$ in $(\eFG)$ by $\tau_1=(0,\tau''_1)C''$ to obtain

\noindent
$(\eFI)$\hfil\hbox{
$\ds\eqalign{
        &\ldots\to\OO^{Z''_n}\buildrel{\tau''_n}\over\to\OO^{Z''_{n-1}}\to\ldots\cr
        &\phantom{12345}\to\OO^{Z''_3}\buildrel{{0}\choose{\tau''_3}}\over\longrightarrow
        \OO^{Z'_1\oplus Z''_2}\buildrel{{\!\!\!1\,\,\,0}\choose{0\,\,\tau''_2}}\over\longrightarrow
        \OO^{Z_1}\buildrel{\tau_1}\over\longrightarrow
        S|U''\to0.\cr
 }
 $}\hfil

	 Both $(\eFH)$ and $(\eFI)$ are \sres{}s by \p\tDD(b).
	 Let $K\subset\OO^{Z_1}$ be the kernel of $\tau_1\in\Hom(\OO^{Z_1},S)(U)$,
	$K'\subset\OO^{Z'_1}$ the kernel of $\tau'_1\in\Hom(\OO^{Z'_1},S)(U')$,
	and $K''\subset\OO^{Z''_1}$ the kernel of $\tau''_1\in\Hom(\OO^{Z''_1},S)(U'')$.
	 Then $K|U'\cong K'\oplus\OO^{Z''_1}$ has an \sres\ by $(\eFH)$, and
	$K|U''\cong\OO^{Z'_1}\oplus K''$ has an \sres\ by $(\eFI)$.
	 Thus $(\eFA)$ is as claimed in part~(a).

	 (b) As the kernel $K$ in (a) satisfies the same conditions as $S$ does,
	we can find by repeated application of (a) \ses{}s
	$0\to K_n\to\OO^{Z_n}\to K_{n-1}\to0$ over $U$ for $n\ge1$,
	where $K_0=S$, and $Z_n$ is a \bspc.
	 Splicing together these \ses{}s we get an \sres\ $(\eDE)$ as claimed.
	 The proof of \t\tFA\ is complete.
}

	 Let $Y$ be a \bspc, $\pi:\CC^N\times\to\CC^N\times Y$ the projection
	$\pi(\za,y)=(\za,0)$, $D\subsub\CC^N$ \pscx\ open, $R:D\to(0,\infty)$
	\cts\ and bounded away from zero with $-\log R$ \psh\ on $D$, and
$$
	\Oa(D,R)=\{(\za,y)\in\CC^N\times Y:\|y\|<R(\za)\}.
\tag\eFJ
 $$
 	 Let $\UU$ be an \scvr\ of $D$, $\Oa=\Oa(D,R)$, and
$$
	\UU(\Oa)=\{U(\Oa)=\pi^{-1}(U)\cap\Oa:U\in\UU\}
\tag\eFK
 $$
 	a \kiemel{basic covering} of $\Oa$.

	 Let $Q=\{\za\in\CC^N:A'_j\le\re\,\za_j\le A''_j,\,B'_j\le\im\,\za_j\le B''_j\}$
	be a compact \kiemel{`cube'} (rectangular box) in $\CC^N$ with $A'_j<A''_j$,
	$B'_j<B''_j$ for $j=1,\ldots, N$.
	 A \kiemel{simple subdivision} of $Q$ into subcubes
	$\{Q_{k_1\ldots k_N}^{l_1\ldots l_N}\}$ is a choice of subdivisions
	$A'_j=a_{j0}<a_{j1}<\ldots<a_{jm}=A''_j$,
	$B'_j=b_{j0}<b_{j1}<\ldots<b_{jm}=B''_j$
	of the edges $[A'_j,A''_j]$, $[B'_j,B''_j]$ of $Q$, $j=1,\ldots,N$,
	where $Q_{k_1\ldots k_N}^{l_1\ldots l_N}=
	\{\za\in\CC^N:a_{j,k_j-1}\le\re\,\za_j\le a_{j,k_j}\,
	b_{j,l_j-1}\le\im\za_j\le b_{j,l_j}\}$ for
	$k_1,\ldots,k_N,l_1,\ldots,l_N=1,\ldots,m$.
	 A \kiemel{simple covering} $\VV=\{V_{k_1\ldots k_N}^{l_1\ldots l_N}\}$
	of $Q$ is an open covering obtained by fattening up the cubes 
	$Q_{k_1\ldots k_N}^{l_1\ldots l_N}$ of a simple subdivision of $Q$ by a
	small amount $0<\epsz<\frac14\min\{a_{j,k}-a_{j,k-1},b_{j,k}-b_{j,k-1}:
	j=1,\ldots,N,\,k=1,\ldots,m\}$, where
	$V_{k_1\ldots k_N}^{l_1\ldots l_N}=\{\za\in\CC^N:a_{j,k_j-1}-\epsz<\re\,\za_j<
	a_{j,k_j}+\epsz,\,b_{j,l_j-1}-\epsz<\im\,\za_j<b_{j,l_j}+\epsz\}$.
	 A \kiemel{\ccvr} $\WW$ of a \pscx\ open subset $\Oa$ of a \bspc\ $X$
	is an \scvr\ $\WW$ of the form
	$\WW=f^{-1}(\VV)=\{f^{-1}(V_{k_1\ldots k_N}^{l_1\ldots l_N})\}$,
	where $f\in\OO(\Oa,\CC^N)$, $Q$ is a cube in $\CC^N$ that contains $f(\Oa)$, and
	$\VV$ is a simple covering of $Q$.
	 (If convenient, we may throw away those 
	$f^{-1}(V_{k_1\ldots k_N}^{l_1\ldots l_N})$ that are empty.)
	 A \ccvr\ plays well with \cpair{}s.

\tetel{\p\tFB.}{With the above notation and hypotheses suppose that $Y$ has a
	\sbs, and \pshdom\ holds in every \pscx\ subset of\/ $\Oa=\Oa(D,R)$,
	and\/ $\UU$ has a (finite) refinement\/ $\WW$ that is a \ccvr\ of $D$.
	 Let $S\to\Oa$ be a sheaf that has an \sres\ over each member of\/ $\UU(\Oa)$.
	 Then $S$ has an \sres\ over\/ $\Oa$.
}

\biz{Proof.}{As this can be proved by the usual induction process of Cousin and
	Cartan relying on \t\tFA\ to amalgamate syzygies over larger and larger
	\cpair{}s, the proof of \p\tFB\ is complete.
}

\tetel{\t\tFC}{Let $X$ be a \bspc\ with a \sbs, $\Oa\subset X$ \pscx\ open,
	$\oa_N\subset\oa_{N+1}\subset\Oa$ \pscx\ open, $N\ge1$, and 
	$\bigcup_{N=1}^\infty\oa_N=\Oa$.
	 Suppose that \pshdom\ holds in every \pscx\ open subset of\/ $\Oa$.
	 Let $S\to\Oa$ be a sheaf \st\ there are \sres{}s
$$
	\ldots\to\OO^{Z^N_n}\buildrel{\tau^{N}_{n}}\over\to\OO^{Z^N_{n-1}}\to\ldots
	\to\OO^{Z^N_1}\buildrel{\tau^{N}_{1}}\over\to S|\oa_N\to0
\tag\eFL
 $$
 	for all $N\ge1$.
	 Then there is a \bspc\ $Z'$ so that the following hold.
\vskip0pt
	{\rm(a)} There is a \ses\/ $0\to K\to\OO^{Z'}\to S\to0$ over\/ $\Oa$
	that is\/ $2$-exact over $\oa_N$ for all $N\ge1$, and $K|\oa_N$ has an
	\sres\ $\ldots\to\OO^{Z'}\to\OO^{Z'}\to K|\oa_N\to0$ with all \bspc{}s
	equal to $Z'$ for all $N\ge1$.
\vskip0pt
	{\rm(b)} There is an \sres\ $\ldots\to\OO^{Z'}\to\OO^{Z'}\to S\to0$ over\/ $\Oa$
	with all \bspc{}s equal to $Z'$.
}

\biz{Proof.}{Let $1\le P<\infty$, and $Z'$ the $\ell_P$-sum of countably infinitely 
	many copies of $Z^N_n$ for $n,N\ge1$.
	 Then $Z'\oplus Z'\cong Z'$, $Z^N_n\oplus Z'\cong Z'$ for $n,N\ge1$, and $\ell_P(Z')\cong Z'$,
	where the isomorphisms are effected by isometries that permute the coordinates.
	 We may easily achieve that in $(\eFL)$ all the \bspc{}s are equal to $Z'$ by
	taking the direct sum of $(\eFL)$ and the following two trivial \sres{}s
$$\eqalign{
	&\ldots\to0\to\OO^{Z'}\buildrel1\over\to\OO^{Z'}\to0\to\OO^{Z'}
	\buildrel1\over\to\OO^{Z'}\to0\to0\cr
	&\ldots\to\OO^{Z'}\buildrel1\over\to\OO^{Z'}\to0
	\to\OO^{Z'}\buildrel1\over\to\OO^{Z'}\to0\to0\to0\cr
} $$
	over $\Oa$.
	 Suppose now that we have \sres{}s
$$
	\ldots\to\OO^{Z'}\buildrel{\tau^N_2}\over\to\OO^{Z'}\buildrel{\tau^N_1}\over\to S|\oa_N\to0
\tag\eFM
 $$
 	with all \bspc{}s equal to $Z'$ for all $N\ge1$.
	 Looking at the image of $(\eFM)$ under the functor $\Hom(\OO^{Z'},-)$ we find,
	by the $2$-exactness of $(\eFM)$ and \t\tDC, for $\tau^{N+1}_1\in\Hom(\OO^{Z'},S)(\oa_N)$
	an $A_N\in\OO(\oa_N,\End(Z'))$ with $\tau^{N+1}_1=\tau^N_1\dot A_N$ over $\oa_N$
	for $N\ge1$.
	 Similarly, we find a $B_N\in\OO(\oa_N,\End(Z'))$ with $\tau^N_1=\tau^{N+1}_1\dot B_N$
	over $\oa_N$ for $N\ge1$.
	 We regard $Z'$ as $\ell_P(Z')$, and similarly as in $(\eFF)$ and $(\eFG)$ by
	adding one or two trivial \sres{}s to $(\eFM)$ we can arrange that the \sres\ $(\eFM)$
	takes the form
$$
	\ldots\to\OO^{Z'}\buildrel{\sa^N_2}\over\to\OO^{Z'}\buildrel{\sa^N_1}\over\to S|\oa_N\to0
\tag\eFN
 $$
 	over $\oa_N$, where the $N$th entry of $\sa^N_1=(0,\ldots,\tau^N_1,0,\ldots)$
	is $\tau^N_1$ for $N\ge1$.

	 Consider the \holo\ operator matrices
	$\tilde A^t_N\in\OO(\oa_N,\GL(\ell_P(Z')))$,
	$\tilde B^t_N\in\OO(\oa_N,\GL(\ell_P(Z')))$
	defined for $0\le t\le 1$ and $x\in\oa_N$ by
$$
	\tilde A^t_N(x)=[\da^j_i+tA_N(x)\da^N_i\da^j_{N+1}]_{i,j\ge1},\quad
	\tilde B^t_N(x)=[\da^j_i-tB_N(x)\da^{N+1}_i\da^j_{N}]_{i,j\ge1},
 $$
 	where $\da^q_p$ is the Kronecker delta.
	 Note that the above operators are in deed invertible since
	$\tilde A^t_N(x)^{-1}=[\da^j_i-tA_N(x)\da^N_i\da^j_{N+1}]_{i,j\ge1}$,
	and $\tilde B^t_N(x)^{-1}=[\da^j_i+tB_N(x)\da^{N+1}_i\da^j_N]_{i,j\ge1}$.
	 Define $\tilde C^t_N\in\OO(\oa_N,\GL(\ell_P(Z')))$ by
	$\tilde C^t_N(x)=\tilde A^t_N(x)\tilde B^t_N(x)$, and note that
	$\sa^N_1\tilde C^1_N=\sa^{N+1}_1$ over $\oa_N$ for $N\ge1$.

	 Indeed,
$$\eqalign{
	[\sa^N_1\tilde C^1_N]_{k\ge1}&=
	\big[\sum_{i,j}\da^i_N\tau^N_1(\da^j_i+A_N\da^N_i\da^j_{N+1})
		(\da^k_j-B_N\da^{N+1}_j\da^k_N)\big]\cr
	&=
	\big[\sum_j\tau^N_1(\da^j_N+A_N\da^j_{N+1})
		       (\da^k_j-B_N\da^{N+1}_j\da^k_N)\big]\cr
	&=\big[\tau^N_1\sum_j(\da^j_N\da^k_j+A_N\cdot\da^j_{N+1}\da^k_j
	                      -B_N\da^k_N\cdot\da^j_N\da^{N+1}_j\cr
	&\phantom{12345678901234567}
			      -A_NB_N\da^k_N\cdot\da^j_{N+1}\da^{N+1}_j)\big]\cr
	&=[\tau^N_1\da^k_N+\tau^{N+1}_1\da^k_{N+1}-0-\tau^N_1\da^k_N]\cr
	&=[\tau^{N+1}\da^k_{N+1}]_{k\ge1}=\sa^{N+1}_1\cr
}
 $$
 	over $\oa_N$, where we used that $\tau^N_1A_NB_N=\tau^N_1$.

	 Let $p\land q=\min\{p,q\}$, and define $\tilde C^t_{pq}\in\OO(\oa_{p\land q},\GL(\ell_P(Z')))$
	for $0\le t\le1$, $x\in\oa_{p\land q}$, and $p,q\ge1$ by
$$
	\tilde C^t_{pq}(x)=
	  \cases
            \tilde C^t_p(x)\tilde C^t_{p+1}(x)\ldots\tilde C^t_{q-1}(x) &\text{if}\quad p<q\cr
            1 &\text{if}\quad p=q\cr
	    (\tilde C^t_q(x)\tilde C^t_{q+1}(x)\ldots\tilde C^t_{p-1}(x))^{-1}&\text{if}\quad p>q
          \endcases.
 $$
 	 Note that $\sa^p_1\tilde C^1_{pq}=\sa^q_1$ over $\oa_{p\land q}$ for $p,q\ge1$.
	 Then $\tilde C^t_{pq}$ is a \holo\ cocycle, i.e.,
	$\tilde C^t_{pq}\tilde C^t_{qr}\tilde C^t_{rp}=1$ over $\oa_{p\land q\land r}$
	for $p,q,r\ge1$ as it is easy to verify.

	 Look at the \holo\ \bvbdl\ $E\to\Oa$ with fiber type $\ell_P(Z')$ whose
	defining cocycle \wrt\ the covering $\{\oa_N:N\ge1\}$ of $\Oa$ is $\tilde C^1_{pq}$.
	 This \bvbdl\ $E$ is topologically trivial over $\Oa$, due to the homotopy
	$\tilde C^t_{pq}$ of cocycles from $\tilde C^1_{pq}$ to $\tilde C^0_{pq}=1$.
	 By \t\tDA(d) our \holo\ \bvbdl{} $E$ is \holo{}ally trivial over $\Oa$, i.e.,
	there are $D_p\in\OO(\oa_p,\GL(\ell_P(Z')))$ \st\
	$\tilde C^1_{pq}(x)=D_p(x)D_q(x)^{-1}$ for $x\in\oa_{p\land q}$, and $p,q\ge1$.
	 Then $\sa^p_1D_p=\sa^q_1D_q$ patch up to a \homo\ $\sa_1\in\Hom(\OO^{\ell_P(Z')},S)(\Oa)$.
	 Replacing $\sa^N_1$ by $\sa_1=\sa^N_1D_N$ in $(\eFN)$ we get another \sres\
$$
	\ldots\to\OO^{Z'}\buildrel{\sa^N_2}\over\to\OO^{Z'}\buildrel{\sa_1}\over\to S|\oa_N\to0
\tag\eFO
 $$
 	of $S|\oa_N$ for $N\ge1$ by \p\tDD(b), where we regard $\ell_P(Z')$ as $Z'$.
	 Let $K=\Ker\,\sa_1\subset\OO^{Z'}$ over $\Oa$, $K_N=\Ker\,\sa^N_1\subset\OO^{Z'}$
	over $\oa_N$, $N\ge1$.
	 Then $K|\oa_N\cong K_N\oplus\OO^{Z'}$ by $(\eFO)$.

	 (b) As the kernel $K$ in (a) satisfies the same conditions as $S$ does, we can find
	by repeated application of (a) \ses{}s $0\to K_n\to\OO^{Z'}\to K_{n-1}\to0$ over $\Oa$
	for $n\ge1$, where $K_0=S$.
	 Splicing together these \ses{}s we get an \sres\ of the type claimed.
	 The proof of \t\tFC\ is complete.

}

\alcim{\sG. GLOBAL (S)-RESOLUTIONS.}

	 In this section we show in many cases that an \ssf\ over a \pscx\ open
	subset of a \bspc\ has a global \sres.

\tetel{\p\tGA.}{Let $X$ be a \bspc\ with a \sbs, $\Oa\subset X$ \pscx\ open,
	$S\to\Oa$ an \ssf, and suppose that \pshdom\ holds in every \pscx\ open
	subset of\/ $\Oa$.
	 Then there is an admissible Hartogs \fn\ $\ba\in\AA$ as in\/ \S\,\sB\ \st\
	there are \sres{}s
$$
	\ldots\to\OO^{Z^N_n}\to\OO^{Z^N_{n-1}}\to\ldots\to\OO^{Z^N_1}\to S|\Oa_N\<\ba\>
\tag\eGA
 $$
 	for all $N\ge1$.
}

\biz{Proof.}{Let $\UU$ be a covering of $\Oa$ by balls $U=B_X(x,r(x))$ for
	$x\in\Oa$ with \cts\ radius \fn\ $r\in C(\Oa,(0,1))$ so small that over
	each $U$ our $S|U$ admits an \sres.
	 By \pshdom\ in $\Oa$ there is an $\aa$ with $10\aa\in\AA$, $10\aa<r$, \st\
	$\BB(\aa)$ as in $(\eBA)$ is a refinement of $\UU$.
	 \p\tBB(f) gives a $\ba\in\AA$, $\ba<\aa$, \st\ the covering 
	$\UU_N=\BB_N(\aa)|\Oa_N\<\ba\>$ has a finite basic refinement $\VV_N$ for all
	$N\ge1$.
	 Any finite basic covering $\VV_N$ of $\Oa_N\<\ba\>$ has a finite refinement
	$\WW_N$ which is a \ccvr\ of $\Oa_N\<\ba\>$ for all $N\ge1$.
	 \p\tFB\ gives us an \sres\ $(\eGA)$ for all $N\ge1$.
	 The proof of \p\tGA\ is complete.
}

	 We now prepare to apply \t\tFC.

\tetel{\p\tGB.}{Let $X$ be a \bspc\ with a \sbs, $\Oa\subset X$ \pscx\ open,
	$\ba\in\AA$, and suppose that \pshdom\ holds in\/ $\Oa$.
	 Then there are $\oa_p\subset\oa_{p+1}\subset\Oa$ \pscx\ open for $p\ge1$
	with\/ $\bigcup_{p=1}\oa_p=\Oa$, and $\oa_p\subset\Oa_p\<\ba\>$ for $p\ge1$.
}

\biz{Proof.}{Define a \fn\ $m:\Oa\to\{1,2,\ldots\}$ as follows.
	 For $x\in\Oa$ let $m(x)$ be the least integer $N\ge1$ \st\
	$x$ has an open \nbd\ $U$ in $\Oa$ with $U\subset\bigcap_{i=0}^\infty\Oa_{N+i}\<\ba\>$.
	 Such a number exists by \p\tBB(a), and $m$ is a locally upper bounded \fn\ of
	$x\in\Oa$ since $m(y)\le m(x)$ if $y\in U$, where $U$ is as above.
	 By \pshdom\ in $\Oa$ we find a \cts\ \psh\ \fn\ $\psi:\Oa\to\RR$
	with $m(x)<\psi(x)$ for $x\in\Oa$.
	 Let $\oa_p=\{x\in\Oa:\psi(x)<p\}$ for $p=1,2,\ldots$.
	 If $x\in\oa_p$, then $m(x)<\psi(x)<p$, so $x\in\Oa_p\<\ba\>$, i.e.,
	$\oa_p\subset\Oa_p\<\ba\>$ for $p\ge1$.
	 As $\oa_p\subset\oa_{p+1}$, and $\bigcup_{p=1}\oa_p=\Oa$
	clearly hold, the proof of \p\tGB\ is complete.
}

\tetel{\t\tGC.}{Let $X$ be a \bspc\ with a \sbs, $\Oa\subset X$ \pscx\ open,
	and $S\to\Oa$ an \ssf.
	 If \pshdom\ holds in every \pscx\ open subset of\/ $\Oa$, then there is
	an \sres\/ $(\eDE)$ over\/ $\Oa$.
}

\biz{Proof.}{\p\tGA\ gives a $\ba\in\AA$ and \sres{}s $(\eGA)$ over $\Oa_N\<\ba\>$
	for $N\ge1$.
	 \p\tGB\ yields $\oa_N\subset\Oa_N\<\ba\>$ for $N\ge1$.
	  Since the restriction of $(\eGA)$ to $\oa_N$ is an \sres\ for $N\ge1$,
	 an application of \t\tFC\ completes the proof of \t\tGC.
}

\alcim{\sH. THE PROOF OF THEOREM~\tAB.}

	 In this section we complete the proof of \t\tAB.

\biz{Proof of \t\tAB.}{Part (a) follows from \t\tEC.

	(b) \p\tEB\ shows that $\OO^E\to M$ is an \ssf\ over $M$,
	the zero extension $(\OO^E)^0\to\Oa$ of which then by \t\tEF\
	is an \ssf\ over $\Oa$.
	 The proof of \t\tAB\ is complete.
}

	 The same reasoning proves the following \t\tHA.

\tetel{\t\tHA.}{With the notation and hypotheses of \t\tAB\ let $S\to M$ be
	an \ssf, and $S^0\to\Oa$ its zero extension to\/ $\Oa$
	 Then $S^0$ is an \ssf\ over\/ $\Oa$, and
	\t\tAC\ holds with\/ $\Oa$ replaced by $M$, and \pscx\ open sets $U$
	replaced by open sets $U$ of the form $U=\tilde U\cap M$, where
	$\tilde U$ is any \pscx\ open subset of\/ $\Oa$.
}

\alcim{\sI. THE PROOF OF THEOREM~\tAC.}

	 In this section we complete the proof of \t\tAC.

\biz{Proof of \t\tAC.}{Part (a) follows from \t\tGC.
	 Part (b) follows from (c) on letting $K=\Ker(\OO^{Z_1}\to S)$.
	 Part (a) follows on applying \t\tDB\ to a global \sres\ in~(c).
	 To prove (d) we see by repeated application of (b) for (S${}'$)-sheaves
	that $S$ admits a global $2$-resolution $(\eAC)$ over $\Oa$,
	hence $S$ is an \ssf\ over $\Oa$.
	 The meaning of (d) is that the class of \ssfs\ is the largest subclass
	of \lcass\ for which \t\tAC(b), i.e., a natural condition, holds.
	 The proof of \t\tAC\ is complete.
}

\alcim{\sJ. THE PROOF OF THEOREM~\tAD.}

	 In this section we complete the proof of \t\tAD.

\biz{Proof of \t\tAD.}{(a) As $I$ is an \ssf\ over $\Oa$ by \t\tAB(a) an
	application of \t\tAC(a) shows the acyclicity of $I$ over $\Oa$.

	(b) As there are local extensions $\tilde f_U$ of $f$, and the
	cocycle $(\tilde f_V-\tilde f_U)$ of $I$ over $\Oa$ can be resolved by (a),
	part (b) is proved.

	(c) Since $(\OO^E)^0$ is an \ssf\ over $\Oa$ by \t\tAB(b),
	we see by \t\tAC(a) that $H^q(\Oa,(\OO^E)^0)$ for $q\ge1$.
	 As $H^q(\Oa,(\OO^E)^0)$ and $H^q(M,\OO^E)$ are canonically
	isomorphic for $q\ge1$, the latter are zero, too.

	(d) As the ideal sheaf $I=I^\CC$ of $M$ in $\Oa$ is an \ssf\ over $\Oa$ by
	\t\tAB(a), there is by \t\tAC(b) a \stes\
$$
	0\to K\to\OO^{Z_1}\buildrel\tau\over\to I\to0
\tag\eJA
 $$
 	over $\Oa$.
	 \p\tCA\ gives a $T\in\OO(\Oa,\Hom(Z_1,\CC))$ that induces the \homo\
	$\tau:\OO^{Z_1}\to\OO$ as $\tau=\dot T$.
	 Let $M'=\{x\in\Oa:T(x)=0\}$.
	 
	 We claim that $M=M'$ as point sets.

	 Let $x_0\in M$ and suppose for a contradiction that $x_0\not\in M'$.
	 As there is a $z_1\in Z_1$ with $T(x_0)z_1\not=0$, there is a small open
	ball $U$ with $x_0\in U\subset\Oa$ \st\ the \fn\ $g\in\OO(U,Z_1)$ defined
	by $g(x)=T(x)z_1$ is nonzero for $x\in U$.
	 As $g\in I(U)$ and $x_0\in M$ we find the contradiction that $g(x_0)=0$.

	 Let $x_0\in M'$ and suppose for a contradiction that $x_0\not\in M$.
	 There is a small open ball $U$ with $x_0\in U\subset\Oa$ that is disjoint
	from the closed set $M$.
	 As the constant $1\in I(U)$, there is a $g\in\OO(U,Z_1)$ with $1=T(x)g(x)$
	for $x\in U$.
	 Letting $x=x_0$ we find the contradiction that $1=0$.

	 Hence $M=M'$, and (d) is proved.

	 The meaning of part (d) is that $M$ can be defined by a global \holo\
	equation $T(x)=0$ in $\Oa$.

	 (e) This follows from (d) upon applying [\rPB, Thm.\,1.2].

	 (f) This follows by the Grauert--Docquier type argument in the proof of
	 [\rPB, Thm.\,6.2]
	 together with (c) and (e).

	  The proof of \t\tAD\ is complete.
}

\alcim{\sK. APPLICATIONS.}

	 In this section we discuss some applications of the theorems in \S\,\sA.
	 See [\rPD, \S\,14] for additional applications.

\tetel{\t\tKA.}{Let $X$ be a \bspc\ with a \sbs, $\Oa\subset X$ \pscx\ open,
	$M\subset\Oa$ a split complex \bsmfd\ of\/ $\Oa$, and $E\to M$ a \holo\
	\bvbdl\ with a \bspc\ $Z$ for fiber type.
	 If \pshdom\ holds in every \pscx\ open subset of\/ $\Oa$, then we
	have the following.
\vskip0pt
	{\rm(a)} Let $Z_1=\ell_p(Z)$, $1\le p<\infty$.
	 Then $E\oplus(M\times Z_1)$ and $M\times Z_1$ are \holo\ isomorphic
	over $M$.
\vskip0pt
	{\rm(b)} $H^q(M,\OO^E)=0$ for $q\ge1$.
\vskip0pt
	{\rm(c)} If $E$ is \cts{}ly trivial over $M$, then
	$E$ is \holo{}ally trivial over $M$.
}

\biz{Proof.}{By \t\tAD(f) there is a \holo\ rectraction $r:\oa\to M$,
	where $\oa$ is \pscx\ open with $M\subset\oa\subset\Oa$.
	 Apply \t\tDA\ to the pull back bundle $r^*E\to\oa$, and then
	restrict back to $M$.
	 That proves (b) and (c).
	 In part (a) we also need to know in advance that $E\oplus(M\times Z_1)$
	is \cts{}ly trivial.  This follows  since
	$r^*(E\oplus(M\times Z_1))=(r^*E)\oplus(\oa\times Z_1)$ 
	is \cts{}ly trivial by [\rPD, Prop.\,7.1].
	 Then (a) follows from this and (c).
	 The proof of \t\tKA\ is complete.
}

	 \p\tKB\ will be useful in the proof of \t\tKD\ below.

\tetel{\p\tKB.}{Let $X$ be a \bspc\ with a \sbs, and $\Oa\subset X$ open.
\vskip0pt
	{\rm(a)} If\/ $\Oa'\subset X$ is \pscx\ open, $\Oa\subset\Oa'$ is open,
	and for each boundary point $x_0$ of\/ $\Oa$ that is not a boundary point
	of\/ $\Oa'$ there is an open set $U\subset X$ with $x_0\in U$, and\/
	$\Oa\cap U$ \pscx\ open in $X$, then\/ $\Oa$ is \pscx\ open in $X$.
\vskip0pt
	{\rm(b)} If\/ $\Oa$ is bi\holo\ to a \pscx\ open subset of $X$, then
	$\Oa$ is \pscx\ open in $X$.
\vskip0pt
	{\rm(c)} If $X=X'\times X''$ is a direct decomposition of \bspc{}s,
	$\pi:X\to X'\times\{0\}$ is the projection $\pi(x',x'')=(x',0)$,
	$\Oa\subset X$ \pscx\ open, $\Oa'\subset X'\times\{0\}$ is \pscx\
	open (relative to $X'$), then $\pi^{-1}(\Oa')\cap\Oa$ is \pscx\ open
	in $X$.
\vskip0pt
	{\rm(d)} If $M\subset\Oa$ is a split complex \bsmfd, and $r:\Oa\to M$ is
	a \holo\ rectraction, then for each $x_0\in M$ there is a ball
	$U=B_X(x_0,\epsz)$ in $X$ and a direct decomposition $X=X'\times X''$
	of \bspc{}s \st\ in $U$ there is a bi\holo{} coordinate system in
	which the rectraction $r$ can be written as a linear projection $\pi$
	as in\/ {\rm(c)}.
\vskip0pt
	{\rm(e)} Let $\Oa\subset X$ be \pscx\ open, $M\subset\Oa$ a split
	complex \bsmfd\ of\/ $\Oa$.
	 Suppose that \pshdom\ holds in every \pscx\ open subset of\/ $\Oa$.
	 Let $D\subset M$ be a relatively open subset of $M$.
	 If in $M$ at every relative boundary point $x_0\in\bd D$
	there is a coordinate ball $U$ in $M$ with $x_0\in U$, and
	$D\cap U$ \cpscx, then there is a \pscx\ open subset $\tilde D$ of\/ $\Oa$
	with $D=\tilde D\cap M$.
\vskip0pt
	{\rm(f)} Let $\Oa\subset X$ be \pscx\ open, $M\subset\Oa$ a split
        complex \bsmfd\ of\/ $\Oa$.
         Suppose that \pshdom\ holds in every \pscx\ open subset of\/ $\Oa$.
         Let $D\subset M$ be a relatively open subset of $M$,
	$E\to D$ be a \cts\ \bvbdl\ with a \bspc\ $Z$ for fiber type, and
	$Z_1=\ell_p(Z)$ for $1\le p<\infty$.
	 Then $E\oplus(D\times Z_1)$ is \cts{}ly isomorphic to $D\times Z_1$.
\vskip0pt
	{\rm(g)} Let $M$ be a complex \bmfd\ modelled on a \bspc\ $X$,
	$Z_1,Z_2$ \bspc{}s, and $g\in\OO(M,\Hom(Z_1,Z_2))$.
	 If $g(x)\in\Hom(Z_1,Z_2)$ is an epimorphism with split kernel for
	each $x\in M$, then the set $K=\Ker\,g=\{(x,\za_1)\in M\times Z_1:
	g(x)\za_1=0\}$ is a \holo\ \bvbdl\ over $M$.
\vskip0pt
	{\rm(h)} Let $X$ be a \bspc, $X^*$ its dual space, $\xi_n\in X^*$,
	$n\ge1$, with\/ $\|\xi_n\|\to\infty$ as $n\to\infty$.
	 Then there is an $x_0\in X$ with\/ $|\xi_n(x_0)|$ unbounded
	as $n\to\infty$.
\vskip0pt
	{\rm(i)} If $U_1,U_2$ are open subsets of a complex \bmfd\ $M$,
	and \pshdom\ holds in $U_1$, and in $U_2$, then \pshdom\ holds
	in $V=U_1\cap U_2$.
\vskip0pt
	{\rm(j)} With the notation and hypotheses of\/ {\rm(e)}
	\pshdom\ holds in $M$.
}

\biz{Proof.}{(a) A well-known criterion for the \pscx{}ity of $\Oa$ in $X$
	runs as follows.
	 An open set $\Oa$ is \pscx\ in $X$ if and only if $X=\Oa$ or else
	for each boundary point $x_0\in\bd\Oa$ in $X$ there are an open set
	$V$ with $x_0\in V$, a sequence of points $x_n\in V\cap\Oa$ with
	$x_n\to x_0$ as $n\to\infty$, and a \holo\ \fn\ $f\in\OO(V\cap\Oa)$
	with $|f(x_n)|$ unbounded as $n\to\infty$.
	 Such a \fn\ $f$ is called a local \holo\ \fn\ on $\Oa$ singular at
	$x_0$.

	 We will use the above criterion to show that our $\Oa$ is \pscx.
	 Let $x_0$ be a boundary point of $\Oa$ in $X$.
	 If $x_0$ is a boundary point of $\Oa'$, then there is a \holo\
	\fn\ $f\in\OO(\Oa')$ that is singular at $x_0$.
	 If $x_0$ is not a boundary point of $\Oa'$, then there is an open set
	$U$ with $x_0\in U$, and $\Oa\cap U$ \pscx\ open in $X$.
	 As $x_0$ is a boundary point of the \pscx\ open set $\Oa\cap U$,
	there is a \holo\ \fn\ $f\in\OO(\Oa\cap U)$ that is singular at $x_0$.
	 Hence $\Oa$ is \pscx\ by the above criterion.

	(b) This is a well-known statement.

	(c) As $\pi^{-1}(\Oa')\cap\Oa$ is the intersection of two \pscx\ open
	subsets $\pi^{-1}(\Oa')=\Oa'\times X''$ and $\Oa$ of $X$, the statement
	follows.

	(d) It is easy to see that at any point $x_0\in M$ the Fr\'echet
	differential $dr(x_0)\in\End(X)$ is a linear projection with split
	kernel.
	 Hence is the statement.

	(e) \t\tAD(f) gives us a \holo\ retraction $r:\Oa'\to M$, where
	$\Oa'$ is \pscx\ open in $X$ with $M\subset\Oa'\subset\Oa$.
	 There is a  \cts\ radius \fn\ $\epsz'\in C(M,(0,1))$ so small that
	over the ball $B_X(x_0,\epsz'(x_0))\subset\Oa'$ the retraction $r$ can
	be linearized to a linear projection by a biholomorphism for all
	$x_0\in M$.
	 By a standard argument with a partition of unitity
	there is a \cts\ radius \fn\ $\epsz''\in C(M,(0,1))$ so small that
	 $\epsz''<\epsz'$, and over the $\epsz''$-balls a doubling inequality
	holds for $\epsz'$, i.e.,
	for all $x_0\in M$ we have that $\epsz''(x_0)<\epsz'(x_0)$, and
	$\frac12\epsz'(z)<\epsz'(y)<2\epsz'(z)$ for all $y,z\in B_X(x_0,\epsz''(x_0))\cap M$.
	 \t\tAD(e) gives a \pscx\ open subset $\oa$ of $X$ with
	$M\subset\oa\subset\{x\in\Oa':\|x-r(x)\|<\frac14\epsz''(r(x))\}$.
	 Let $\tilde D=r^{-1}(D)\cap\oa$.
	 Then $\tilde D$ is an open subset of $\Oa$, and
	$\tilde D\cap M=D$ since if $x\in\tilde D\cap M$, then
	$r(x)=x\in D$, and if $x\in D$, then 
	$x\in r^{-1}(D)\cap M\subset\tilde D\cap M$.

	 We proceed to show via (a) that $\tilde D$ is \pscx\ open in $X$.
	 If $\tilde D=X$, then $\tilde D$ is \pscx\ open in $X$.
	 If $\tilde D\not=X$, then let $x_0$ be any boundary point of
	$\tilde D$.
	 As $\tilde D$ is an open subset of the \pscx\ open subset $\Oa'$
	of $X$ in order to apply (a) it is enough to check that at any boundary
	point $x_0$ of $\tilde D$ in $X$ that is not a boundary point of
	$\Oa'$ there is an open set $U\subset X$ \st\ $\tilde D\cap U$ is \pscx\
	open in $X$, and $x_0$ is a limit point of $\tilde D\cap U$, but not
	a point of $\tilde D\cap U$.
	 As $x_0$ is not a boundary point of $\Oa'$, i.e., $x_0\in\Oa'$,
	our retraction $r$ is defined at $x_0$.
	 Since there are points $x_n\in\tilde D$ with $x_n\to x_0$ as $n\to\infty$,
	we see that $r(x_n)\to r(x_0)$ as $n\to\infty$, i.e., $r(x_0)$ is
	in the closure of $D$ relative to $M$.

	 We claim that if $V\subset B_X(r(x_0),\frac14\epsz''(r(x_0)))\cap M$, then
	$r^{-1}(V)\cap\oa\subset B_X(r(x_0),\epsz'(r(x_0)))$.

 	 Indeed, we must show for $x\in r^{-1}(V)\cap\oa$ that $\|x-r(x_0)\|<\epsz'(r(x_0))$.
	 As $r(x)\in V$, and $x\in\oa$, we have the inequalities
$$
	\|r(x)-r(x_0)\|<{\ts\frac14}\epsz''(r(x_0)),\qquad
	\|x-r(x)\|<{\ts\frac14}\epsz''(r(x)),
 $$
 	adding up which implies that
$$\eqalign{
	\|x-r(x_0)\|&<{\ts\frac14}\epsz''(r(x_0))+{\ts\frac14}\epsz''(r(x))\cr
	&<{\ts\frac14}\epsz'(r(x_0))+{\ts\frac14}\epsz'(r(x))\cr
	&<{\ts\frac14}\epsz'(r(x_0))+{\ts\frac24}\epsz'(r(x_0))\cr
	&<\epsz'(r(x_0)),\cr
}
 $$
 	where we applied in the penultimate inequality the doubling
	property of $\epsz'$ on the ball $B_X(r(x_0),\frac14\epsz''(r(x_0)))\cap M$.

	 If $r(x_0)\in D$, then the point $r(x_0)$ is contained
	in a coordinate ball
	$V\subset B_X(r(x_0),\frac14\epsz''(r(x_0)))\cap D\subset D$ relative to $M$.
	By the claim above the set $U=r^{-1}(V)\cap\oa$ 
	is contained in a ball $B_X(r(x_0),\epsz'(r(x_0)))$
	in which $r$ can be linearized to a linear projection.
	 Then $U=\tilde D\cap U$ is \pscx\ open in $X$ by
	(b) and (c).
	 We now show that $x_0$ is limit point of $\tilde D\cap U$.
	 Indeed, 
	 let $x_n\in\tilde D$ be any sequence with $x_n\to x_0$ as $n\to\infty$.
	 As $r(x_n)\to r(x_0)\in V$, and $V$ is open relative to $M$,
	there is an $N$ with $r(x_n)\in V$ for all $n\ge N$.
	 So $x_n\in r^{-1}(V)\cap\oa=U$ for $n\ge N$, i.e.,
	$x_0$ is a limit point of $U=\tilde D\cap U$.

	 If $r(x_0)\in\bd D$ is in the boundary of $D$ relative to $M$,
	then there is a coordinate ball $V\subset B_X(r(x_0),\epsz''(r(x_0)))\cap M$
	relative to $M$ with $x_0\in V$ and $V\cap D$ \cpscx\ open relative to $M$.
	By the claim above the set $U=r^{-1}(V)\cap\oa$ is
	contained in a ball $B_X(r(x_0),\epsz'(r(x_0)))$ in
	which $r$ can be linearized to a linear projection.
	 Then $U$ and $\tilde D\cap U=r^{-1}(V\cap D)\cap\oa$ are \pscx\ open in $X$
	by (b) and (c).
	 We now show that $x_0$ is a limit point of $\tilde D\cap U$.
	 Let $x_n\in\tilde D=r^{-1}\cap\oa$ be any sequence with
	$x_n\to x_0$ as $n\to\infty$.
	 As $r(x_n)\to r(x_0)\in D\cap V$, and $D\cap V$ is open relative to $M$,
	there is an $N$ with $r(x_n)\in D\cap V$ for all $n\ge N$.
	 Then $x_n\in r^{-1}(D\cap V)\cap\oa=\tilde D\cap U$
	for $n\ge N$, i.e., $x_0$ is a limit point of $\tilde D\cap U$.
	 Thus as our $\tilde D$ is \pscx\ open in $X$ by (a), the proof of (e)
	is complete.

	(f) By \t\tAD(f) there is a \holo\ retraction $r:\oa\to M$, where
	$\oa$ is \pscx\ open $X$ with $M\subset\oa\subset\Oa$.
	 Look at $r^*(E\oplus(D\times Z_1))=(r^*E)\oplus(r^{-1}(D)\times Z_1)$
	and apply [\rPD, Prop.\,7.1], then restrict back to $D$.

	(g) This follows from the inverse function theorem for \holo\
	maps of \bspc{}s.  Note that any closed \fdml\ or finite codimensional
	subspace of $Z_1$ is split, and so is any closed subspace of $Z_1$ 
	if $Z_1$ is a \hspc.

	(h) This follows from the principle of uniform boundedness or
	the principle of condensation of singularities in linear functional analysis.

	(i) Let $u:V\to\RR$ be a locally upper bounded \fn, $U=U_1\cup U_2$,
	and $\chi_1,\chi_2:U\to\RR$ a \cts\ partition of unity subordinate to
	the open covering $\{U_1,U_2\}$ of $U$.
	 As $u\chi_i$ is a locally upper bounded \fn\ on $U_i$, \pshdom\ in $U_i$
	gives a \cts\ \psh\ \fn\ $\psi_i:U_i\to\RR$ with $u(x)\chi_i(x)<\psi_i(x)$
	for $x\in U_i$, $i=1,2$.
	 As $u(x)=u(x)\chi_1(x)+u(x)\chi_2(x)<\psi_1(x)+\psi_2(x)$ for $x\in V$,
	the \cts\ \psh\ \fn\ $\psi_1+\psi_2:V\to\RR$ dominates $u$ on $V$.

	(j) Let $u:M\to\RR$ be a locally upper bounded \fn, and define $u':\Oa\to\RR$
	by $u'(x)=u(x)$ if $x\in M$, and $u'(x)=0$ if $x\in\Oa\setminus M$.
	 As $u'$ is easily seen locally upper bounded (since $M$ is relatively closed
	in $\Oa$), \pshdom\ in $\Oa$ gives a \cts\ \psh\ \fn\ $\psi':\Oa\to\RR$
	with $u'(x)<\psi'(x)$ for $x\in\Oa$.
	 Then $\psi=\psi'|M$ is a \cts\ \psh\ \fn\ on $M$ that dominates $u$.

	 The proof of \p\tKB\ is complete.
}

	 Among \fdml\ complex \mfd{}s the class of Stein \mfd{}s can be
	characterized by cohomological criteria.
	 There are also cohomological criteria for open subsets of a Stein
	\mfd\ $M$ to be themselves Stein.
	 Here is one such criterion by Leiterer.

\tetel{\t\tKC.}{{\rm(Leiterer, [\rLtB])}
	Let $M$ be a Stein \mfd\ of complex dimension $n$, and $D\subset M$ open.
	Then the following are equivalent.
\vskip0pt
	{\rm(a)} $D$ is a Stein \mfd.
\vskip0pt
	{\rm(b)} $H^1(D,\OO)=0$, and any topologically trivial \holo\ \vbdl\
	over $D$ is \holo{}ally trivial over $D$.
\vskip0pt
	{\rm(c)} $H^1(D,\OO)=0$, and for every corank\/ $1$ \holo\ vector subbundle
	$E$ of $D\times\CC^{2n+1}$ \st\ for some $m$ the bundle
	$E\oplus(D\times\CC^m)$ is topologially trivial over $D$,
	there is a topologically trivial \holo\
	\vbdl\ $F\to D$ with $E\oplus F$ \holo{}ally trivial over $D$.
\vskip0pt
	{\rm(d)} $H^1(D,\OO^E)=0$ for every corank\/ $1$ \holo\ vector subbundle
	$E$ of $D\times\CC^{2n+1}$.
\vskip0pt
	{\rm(e)} For every choice of \holo\ \fns\ $g_1,\ldots,g_{2n+1}\in\OO(D)$
	without common zeros in $D$ there are \holo\ \fns\ $f_1,\ldots,f_{2n+1}\in\OO(D)$
	with $\sum_{i=1}^{2n+1}f_i(x)g_i(x)=1$ for $x\in D$.
}

	 We give in \t\tKD\ below an analog of \t\tKC\ above.

\tetel{\t\tKD.}{Let $X$ be a \bspc\ with a \sbs, $\Oa\subset X$ \pscx\ open,
	$M$ a split complex \bsmfd\ of\/ $\Oa$, and $D\subset M$ relatively open.
	 Suppose that \pshdom\ holds in every \pscx\ open subset of\/ $\Oa$.
	 Then the following are equivalent.
\vskip0pt
	{\rm(a)} There is a \pscx\ open subset $\tilde D$ of $X$ with
		$D=M\cap\tilde D$.
\vskip0pt
	{\rm(b)} $H^1(D,\OO^Z)=0$ for any \bspc\ $Z$, and any \cts{}ly
	trivial \holo\ \bvbdl\ over $D$ is \holo{}ally trivial over $D$.
\vskip0pt
	{\rm(c)} $H^1(D,\OO^Z)=0$ for any \bspc\ $Z$, and for any corank\/ $1$
	\holo\ Banach vector subbundle $E$ of $D\times X$ over $D$ there is a
	\holo\ \bvbdl\ $F\to D$ with $E\oplus F$ \holo{}ally trivial over $D$.
\vskip0pt
	{\rm(d)} $H^1(D,\OO^E)=0$ for every corank\/ $1$ Banach vector subbundle
	$E$ of $D\times X$ over $D$.
\vskip0pt
	{\rm(e)} For every $g\in\OO(D,X)$ with $g(x)\not=0$ for $x\in D$ there is
	an $f\in\OO(D,X^*)$ with $f(x)\cdot g(x)=1$ for $x\in D$, where $X^*$ is
	the \bspc\ dual to $X$, and the dot denotes the natural pairing
	$X^*\times X\to\CC$.
\vskip0pt
	{\rm(f)} Plurisubharmonic domination holds in $D$.
}

\biz{Proof.}{(a${}\Rightarrow{}$b) As $D=M\cap\tilde D$ we see that $D$ is a split
	complex \bsmfd\ of the \pscx\ open set $\tilde D$ in $X$.
	 Thus (b) follows from \t\tKA(b).

	 (a${}\Rightarrow{}$c) Let $F=D\times Z_1$.
	 Then (c) holds by \t\tKA(a).

	 (a${}\Rightarrow{}$d) Part (d) follows from \t\tKA(b).

	 (a${}\Rightarrow{}$e) As $E=\Ker\,g=\{(x,\xi)\in D\times X^*:\xi g(x)=0\}$
	 is the kernel of the epimorphism (with split kernel) in $\OO(D,\Hom(X^*,\CC))$
	 defined by $(x,\xi)\mapsto\xi g(x)$, \p\tKB(g) shows that
	 $E$ is a \holo\ \bvbdl\ over $D$.
	  As $g(x_0)\not=0$ for $x_0\in D$, there is a \cts\ linear functional
	 $\xi_0\in X^*$ with $\xi_0 g(x_0)=1$.
	  Then there is an open \nbd\ $U_{x_0}$ of $x_0$ in $D$ with $\xi_0 g(x)\not=0$ 
	 for $x\in U_{x_0}$.
	  Define $f_{x_0}\in\OO(U,X^*)$ by $f_{x_0}(x)=\frac{1}{\xi_0 g(x)}\xi_0$.
	  Then $f_{x_0}\cdot g=1$ on $U_{x_0}$.
	  Let $\UU=\{U_{x_0}:x_0\in D\}$, and look at the cocycle 
	 $(f_{y_0}-f_{x_0})\in Z^1(\UU,\OO^E)$.
	  As $H^1(D,\OO^E)=0$ by \t\tKA(b) we have a cochain $h_{x_0}\in C^0(\UU,\OO^E)$
	 with $f_{y_0}-f_{x_0}=h_{y_0}-h_{x_0}$ on $U_{x_0}\cap U_{y_0}$.
	  Then $f=f_{x_0}-h_{x_0}=f_{y_0}-h_{y_0}$ patch up to a well defined \fn\
	 $f\in\OO(D,X^*)$ with $f(x)g(x)=f_{x_0}(x)g(x)-h_{x_0}(x)g(x)=1-0$ for $x\in U_{x_0}$,
	 $x_0\in D$.
	  Thus (e) follows.

	 (a${}\Rightarrow{}$f) As $D=\tilde D\cap M$ \p\tKB(j) shows that
	 \pshdom\ holds in $D$.

	  (b${}\Rightarrow{}$c) Let $F=D\times Z_1$, and apply \p\tKB(f).
	  Then $E\oplus F\to D$ is \cts{}ly trivial, hence it is \holo{}ally
	 trivial by (b), and so (c) follows.

	  (d${}\Rightarrow{}$e) See the proof of (a${}\Rightarrow{}$e) above.

	  (e${}\Rightarrow{}$a) At any boundary point $x_0$ of $D$ relative to $M$ there is
	  a \holo\ \fn\ $\fii\in\OO(D)$ that is singular at $x_0$.
	   Indeed, look at $g\in\OO(D,X)$ defined by $g(x)=x-x_0$.
	   As $g(x)\not=0$ for $x\in D$, there is an $f\in\OO(D,X^*)$ with
	  $f(x)g(x)=1$ for $x\in D$.
	   Let $x_n\in D$, $n\ge1$, be any sequence of points with $x_n\to x_0$ as
	  $n\to\infty$.
	   As $g(x_n)\to0$ as $n\to\infty$, we find that $\|f(x_n)\|$ may not be
	  bounded as $n\to\infty$ (since otherwise $1=f(x_n)g(x_n)\to0$ as $n\to\infty$
	  would hold).
	   \p\tKB(h) gives a point $\tilde x\in X$ with $|f(x_n)\tilde x|$ unbounded as
	  $n\to\infty$.
	   Define $\fii\in\OO(D)$ by $\fii(x)=f(x)\tilde x$.
	   Then $\fii$ is singular at $x_0$.

	   Let $U$ be a coordinate ball relative to $M$ with $x_0\in U$.
	   Then $U\cap D$ is \cpscx\ open relative to $M$.
	   Indeed, any boundary point $y_0$ of $U\cap D$ relative to $M$
	  is a boundary point of $U$ relative to $M$ or a boundary point of
	  $D$ relative to $M$.
	   In either case there is a \holo\ \fn\ $\fii\in\OO(U)$ or
	  $\fii\in\OO(D)$ that is singular at $y_0$.
	   Thus $U\cap D$ is \cpscx\ open relative to $M$, as claimed,
	  by the criterion in the proof of \p\tKB(a).

	   \p\tKB(e) thus applies and gives us a \pscx\ open $\tilde D$ in $X$
	  with $D=M\cap\tilde D$, completing the proof of (e${}\Rightarrow{}$a).

	(c${}\Rightarrow{}$a) The bundles $E=\Ker\,g$ introduced in the proof
	of (a${}\Rightarrow{}$e) are corank $1$ Banach vector subbundles of $D\times X$,
	and (c) provides a \holo\ \bvbdl\ $F\to D$ with $E\oplus F\cong D\times Z$
	\holo{}ally trivial over $D$, where $Z$ is a \bspc.
	 Thus by (c) we see that $0=H^1(D,\OO^Z)=H^1(D,\OO^E)\oplus H^1(D,\OO^F)$,
	i.e., $H^1(D,\OO^E)=0$.
	 Thus (c${}\Rightarrow{}$e), and as (e${}\Rightarrow{}$a), we find that
	(c${}\Rightarrow{}$a).

	(f${}\Rightarrow{}$a) Let $x_0\in\bd D$ be any boundary point of $D$ relative
	to $M$.
	 There is a small enough ball $U=B_X(x_0,r)$ in $\Oa$ with $U\cap M$
	\cpscx\ open relative to $M$.
	 Then \pshdom\ holds in $U\cap M$ by \p\tKB(j), and in $(U\cap M)\cap D=U\cap D$
	by \p\tKB(i), hence $U\cap D$ is \cpscx\ open relative to $M$.
	 Thus \p\tKB(e) applies and gives a \pscx\ open subset $\tilde D$ of $\Oa$ with
	$D=\tilde D\cap M$.

	 The proof of \t\tKD\ is complete.
}

\tetel{\t\tKE.}{Let $M$ be a complex \bmfd\ modelled on a \bspc\ $X$ with a \sbs,
	$S\to M$ an \ssf, and\/ $\UU$ an \scvr\ of $M$.
	 Suppose that \pshdom\ holds in every \pscx\ open subset of $X$.
	 Then the following hold.
\vskip0pt
	{\rm(a)} The covering\/ $\UU$ is a Leray covering of $M$ for the
	sheaf $S$, and $H^q(M,S)$ is naturally isomorphic to the cohomology
	group $H^q(\UU,S)$ of alternating cochains of\/ $\UU$ for all $q\ge1$.
\vskip0pt
	{\rm(b)} If\/ $\UU$ is finite, say, it has $n=1,2,\ldots$ elements,
	then $H^q(M,S)=0$ for $q\ge n$.
\vskip0pt
	{\rm(c)} Let $M\subset X$ be \pscx\ open, $f_i\in\OO(M)$,
	$U_i=\{x\in M:f_i(x)\not=0\}$, $i=1,\ldots,n$, $n\ge1$,
	$\UU=\{U_i:i=1,\ldots,n\}$, $M_0=\bigcap_{i=1}^n(M\setminus U_i)$,
	and $S_0\to M\setminus M_0$ an \ssf.
	 Then $H^q(M\setminus M_0,S_0)$ is naturally isomorphic to
	$H^q(\UU,S_0)$ for $q\ge1$, and is, in particular, zero for $q\ge n$.
}

\biz{Proof.}{As (a) follows from \t\tAC(a), while (b) from (a), and (c)
	from (b) on looking at the \scvr\ $\UU$, the proof of \t\tKE\ is complete.
}

\alcim{\sL. DISCUSSION.}

	 In this section we make some informal remarks on the methods adopted in
	this paper.

	 To prove cohomology vanishing for analytic sheaves via the method of
	amalgamation of syzygies it is necessary to be able to obtain local
	resolutions, to paste local resolutions to make resolutions over bigger
	sets, and to pass to the limit.
	 The Oka coherence theorem is the main tool in \fdms\ to obtain local
	resolutions.
	 No analog of that fundamental theorem of Oka seems to be known currently
	in \idms.
	 To obtain our local resolutions we simply define away the problem.
	 To paste local resolutions together we need to lift local epimorphisms
	to \holo\ operator valued \fns\ through another local epimorphism.
	 Once again we define away the difficulty in the notion of $2$-exactness.
	 Once the liftings are in hand we need the fact that a \holo\ \bvbdl\
	that is topologically trivial is also \holo{}ally trivial in
	certain cases.
	 This was one of the main raison d'\^etre of [\rPD], where this
	was shown.

	 To pass to the limit there are two basic ways.
	 One requires us to show that our sheaf is acyclic over the members
	of the exhaustion and to pass to the limit in the first cohomology of
	the sheaf at hand.
	 This could be done in either of two ways.
	 The one is to do cohomology with bounds for the sheaf --- this seems
	very difficult in \idms.
	 The other is to use the multiplicative Runge type approximation
	Hypothesis$(X,\GL(Z))$ in [\rPA], where $X,Z$ are \bspc{}s, which
	hypothesis is still not proved.
	 The second way of passage to the limit is to obtain a global
	resolution.
	 This is the path that we took here.
	 The reason for that is that we can reduce the problem to proving that
	a topologically trivial \holo\ \bvbdl\ is \holo{}ally trivial in certain
	cases, i.e., to something already available.
	 Once a global resolution is available the proof of the acyclicity of
	the sheaf can be effected through the use of the ensuing dimension
	shifting formula, which in \fdms\ is sufficient by itself, while in
	\idms\ it can be combined with exhaustions whose members have finite Leray
	coverings.
	 This was done in [\rPC] explicitly for this purpose.

	 The notion of $2$-exactness could be fine tuned by requiring
	exactness under the functor $\Hom(\OO^Z,-)$ for not all \bspc{}s $Z$,
	but just for a class of them, which suffices for local resolutions,
	and local amalgamations, and in the resolutions we could allow only
	special type of epimorphisms.
	 E.g., if $S\to\Oa$ is an analytic sheaf over a \pscx\ open set $\Oa$
	that admits local $1$-resolutions $(\eDE)$ with all \bspc{}s $Z_n$
	\fdml, then we could look at $\Hom(\OO^Z,-)$ with $Z$ \fdml\ only,
	in which case the lifting property is trivial.
	 In the limit we could allow ourselves to wind up with a slightly
	different type of global resolution, e.g., in the above example
	we can produce a global $1$-resolution
$$
	\ldots\to\OO^{\ell_2}\buildrel{\tau_2}\over\to\OO^{\ell_2}\buildrel{\tau_1}\over\to S\to0
 $$
 	over $\Oa$, where, locally or over $U=\oa_p$ as in \t\tGC,
	each map $\tau_n$ is of the form
	$\tau_n=\diag(0,\sa_n)A_n:\OO^{\ell_2\oplus\CC^{N_n}}\to\OO^{\ell_2\oplus\CC^{N_{n-1}}}$,
	where $A_n\in\OO(U,\GL(\ell_2\oplus\CC^{N_n}))$,
	$\sa_n\in\Hom(\OO^{N_n},\OO^{N_{n-1}})(U)$ for $n\ge2$,
	and $\tau_1=(0,\sa_1)A_1:\OO^{\ell_2\oplus\CC^{N_n}}\to S$,
	$\sa_1\in\Hom(\OO^{N_1},S)(U)$, and $N_n$ are nonnegative integers for $n\ge1$.

	 Some words on analytic subsets of \bspc{}s are in order.
	 The notion of a Banach analytic set in the most general
	sense seems pathological.
	 Douady [\rD] showed that any compact metric space can be embedded
	in a suitable \bspc\ as a Banach analytic set.
	 Later Pestov [\rPsA, \rPsB] proved the same for any complete metric space.
	 Even closed linear subspaces may raise problems.
	 Let $A$ be a closed linear subspace of a \bspc\ $X$.
	 Is $A$ a Banach analytic subset of $X$?
	 Formally, of course, yes, since one can write
	$A=\{\pi=0\}$, where $\pi:X\to X/A$ is the natural projection.
	 If the associated \ses\ $0\to A\to X\to X/A\to0$ splits,
	then well and good, $A$ has a direct complement.
	 If the above sequence does not split, can one still understand
	some of the complex analytic properties of $A$ from those of $X$?
	 E.g., can we understand much about any separable \bspc\ $A$ by
	embedding it into a universal space like $X=C[0,1]$?
	 Well, not most directly, but perhaps \pshdom\ in $C[0,1]$,
	if available, would have a positive impact on the study of $A$.
	 A nice class of Banach analytic sets seems to be that of
	the split complex \bsmfd{}s, (and possibly their finite branched coverings),
	which can be well studied, as demonstrated in this paper, with
	current technology.

	 The results of this paper naturally raise some questions
	about \ssfs\ and complex \bmfd{}s.
	 Here is one for instance.
	 Let $D$ be a complex \hmfd, and suppose that $H^q(D,S)=0$
	for any $q\ge1$ and any \ssf\ $S$ over $D$.
	 Can $D$ be \holo{}ally embedded as a closed complex \hsmfd\ in the \hspc\ $\ell_2$?
	 Recently, Aaron B.~Zerhusen has shown, based on the notion of
	\holo\ domination in [\rLB], that if $D$ is a \pscx\ open
	subset of $\ell_2$, then $D$ can be \holo{}ally embedded in $\ell_2$
	as a \hsmfd.
	 So in particular, if $D$ is as in \t\tKD\ with $X=\ell_2$, then
	the answer to the above question for $D$ is yes.

	 In conclusion we would like to remark that parts of the present
	paper, especially Theorems~\tFA\ and \tFC,
	were inspired by [\rLtA] by Leiterer.
	 While we could not manage to cite any part of it, we are
	glad that we could at least quote [\rLtB] by him instead.
	
\comment	
	{\vekony Disclosure.}
	 Upon reading the paper [Lempert-Vbdl] by Lempert in March 2003 the author immediately
	realized that Lempert's exhaustion method would allow for a 
	vanishing theorem for a fairly general class of analytic sheaves.
	 During the summer of 2003 the author received employment as an assistant
	visiting researcher at the Department of Mathematics at the University
	of California at Riverside in part for developing such a general
	vanishing theorem.
	 The period from April 2003 to November 2004 was a very difficult
	one for the author during which while he could think about the
	problems of vanishing, and could develop the main ideas in notes,
	he could not fight off a strong feeling of futility, and find the
	necessary strength of heart to produce the required papers.
	 They were produced in December 2004 through July 2005, and
	comprise [Koszul, Vbdl], and the present paper.
	 In October 2004 Lempert, the author's advisor for 1996--2000,
	announced at a conference in New Mexico that he had obtained a
	general sheaf cohomology vanishing theorem.
	 The author has no knowledge of the statement or proof of that
	theorem of Lempert, nor has he received any mathematical
	communication from Lempert on the matter.
	 Given the promise that he made and the support that he received,
	the author, nevertheless, felt duty bound to develop his ideas,
	no matter how belated.
	 He is grateful to Professor L\'aszl\'o Lempert for delaying
	the submission of his manuscript on the topic until the end of
	July 2005 to allow the author time to finish the present manuscript.
\endcomment

	{\vekony Acknowledgements.} The work on this paper started at the Riverside
	campus of the University of California, continued at its San Diego
	campus, and was completed at Georgia State University.
	 The author is grateful to these institutions and to 
	Professors Bun Wong of UCR, Peter Ebenfelt of
	UCSD, and Mih\'aly Bakonyi of GSU, whose support has been essential to him.

\vskip0.30truein
\centerline{\scVIII References}
\vskip0.20truein
\baselineskip=11pt
\parskip=7pt
\frenchspacing
{\rmVIII

	[\rD] Douady, A.,
	{\itVIII A remark on Banach analytic spaces},
	Symposium on infinite-di\-men\-sion\-al topology
	(Baton Rouge, LA, 1967),
	41--42,
	Ann. of Math. Stud., 69,
	Princeton Univ. Press,
	Princeton, NJ, 1972.

	[\rG] Gunning,~R.C.,
	{\itVIII
	Introduction to holomorphic functions of several variables},
	Vol. III,
	Wadsworth \& Brooks/Cole, Belmont, California, 
	(1990).

	[\rLtA] Leiterer,~J.,
	{\itVIII
	Banach coherent analytic Fr\'echet sheaves},
	Math.{} Nachr., {\bfVIII 85} (1978), 91--109.

	[\rLtB] \vonal,
	{\itVIII
	Equivalence of Steinness and validity of Oka's principle
	for subdomains of Stein manifolds},
	Math.{} Nachr., 
	{\bfVIII 89}
	(1979),
	181--183.

	[\rLA] Lempert,~L.,
    	{\itVIII
   	The Dolbeault complex in infinite dimensions~III},
    	Invent.{} Math., {\bfVIII 142} (2000), 579--603.

	[\rLB] \vonal,
	{\itVIII
	Plurisubharmonic domination},
	J. Amer. Math. Soc.,
	{\bfVIII 17}
	(2004),
	361--372.

	[\rLC] \vonal,
	{\itVIII
	Vanishing cohomology for holomorphic vector bundles in a Banach setting},
	Asian J. Math.,
	{\bfVIII 8}
	(2004),
	65--85.

    	[\rLD] \vonal,
    	{\itVIII
	Acyclic sheaves in Banach spaces},
	Contemporary Math.,
	{\bfVIII 368}
	(2005),
	313--320.

	[\rN] Noverraz,~P.,
	{\itVIII
	Pseudo-convexit\'e, convexit\'e polynomiale et domains
	d'holomorphie en dimension infinie},
	North--Holland, Amsterdam, (1973).
	
	[\rPA] Patyi,~I.,
	{\itVIII
	On the Oka principle in a Banach space I},
	Math. Ann.,
	{\bfVIII 326}
	(2003),
	417--441.

	[\rPB] \vonal,
	{\itVIII
	Analytic cohomology of complete intersections in a Banach space},
	Ann. Inst. Fourier (Grenoble),
	{\bfVIII 54}
	(2004),
	147--158.

	[\rPC] \vonal,
	{\itVIII
	An analytic Koszul complex in a Banach space},
	manuscript, (2005).

	[\rPD] \vonal,
	{\itVIII
	On holomorhic Banach vector bundles over Banach spaces},
	manuscript, (2005).

	[\rPsA] Pestov, V.,
	{\itVIII Douady's conjecture on Banach analytic spaces},
	C. R. Acad. Sci. Paris S\'er. I Math,
	{\bfVIII 319}
	(1994),
	no. 10,
	1043--1048.

	[\rPsB] \vonal,
	{\itVIII Analytic subsets of Hilbert spaces},
	Colloque Trajectorien \`a la M\'emoire de George Reeb et
	Jean-Louis Callot (Strasbourg-Obernai, 1995),
	75--80,
	Pr\'epubl. Inst. Rech. Math. Av., 1995/13,
	Univ. Louis Pasteur, Strasbourg, 1995.
}
\vskip0.20truein
\centerline{\vastag*~***~*}
\vskip0.15truein
{\scVIII
	Imre Patyi,
	Department of Mathematics and Statistics,
	Georgia State University,
	Atlanta, GA 30303-3083, USA,
	{\ttVIII ipatyi\@gsu.edu}}
\bye